\documentclass{amsart}

\usepackage[english]{babel}

\usepackage[a4paper,top=2cm,bottom=2cm,left=3cm,right=3cm,marginparwidth=1.75cm]{geometry}

\usepackage{amsmath}
\usepackage{graphicx}
\usepackage{amsfonts}
\usepackage{amsthm}
\usepackage{caption}

\usepackage[colorlinks=true, allcolors=blue]{hyperref}
\newtheorem{theorem}{Theorem}[section]
\newtheorem{proposition}[theorem]{Proposition}

\newtheorem{lemma}[theorem]{Lemma}
\theoremstyle{remark}
\newtheorem{remark}[theorem]{Remark}
\theoremstyle{definition}
\newtheorem{definition}{Definition}[section]

\DeclareMathOperator{\E}{\mathbb{E}}

\DeclareMathOperator{\N}{\mathbb{N}}
\DeclareMathOperator{\R}{\mathbb{R}}
\DeclareMathOperator{\cF}{\mathcal{F}}

\DeclareMathOperator{\cT}{\mathcal{T}}

\DeclareMathOperator{\cN}{\mathcal{N}}

\DeclareMathOperator{\cD}{\mathcal{D}}
\DeclareMathOperator{\cB}{\mathcal{B}}
\DeclareMathOperator{\cS}{\mathcal{S}}

\DeclareMathOperator{\fT}{\mathfrak{T}}
\DeclareMathOperator{\bP}{\mathbb{P}}

\newcommand{\der}[2]{\frac{d #1}{d #2}}

\newcommand{\Norm}[2]{\left\Vert #1 \right\Vert_{#2}}

\title{A sojourn-based approach to semi-Markov Reinforcement Learning}
\author{Giacomo Ascione}
\address{Scuola Superiore Meridionale}
\email{giacomo.ascione@unina.it}
\author{Salvatore Cuomo}
\address{Dipartimento di Matematica e Applicazioni, Universit\'{a} degli Studi di Napoli Federico II}
\email{salvatore.cuomo@unina.it}
\begin{document}
\maketitle

\begin{abstract}
In this paper we introduce a new approach to discrete-time semi-Markov decision processes based on the sojourn time process. Different characterizations of discrete-time semi-Markov processes are exploited and decision processes are constructed by their means. With this new approach, the agent is allowed to consider different actions depending also on the sojourn time of the process in the current state. A numerical method based on $Q$-learning algorithms for finite horizon reinforcement learning and stochastic recursive relations is investigated. Finally, we consider two toy examples: one in which the reward depends on the sojourn-time, according to the \textit{gambler's fallacy}; the other in which the environment is semi-Markov even if the reward function does not depend on the sojourn time. These are used to carry on some numerical evaluations on the previously presented $Q$-learning algorithm and on a different naive method based on deep reinforcement learning.
\end{abstract}

\section{Introduction}
Reinforcement Learning plays a prominent role in Machine Learning nowadays. In this paradigm, the intelligent agent decides the action that it has to take (which could influence or not the environment) in order to maximize a cumulative reward (or minimize a cumulative cost). Practically, it consists in the adaptation of the notion of \textit{Operant conditioning} (see, for instance, \cite{staddon2003operant}) to artificial intelligence.\\
Reinforcement Learning is usually formalized by means of Markov Decision Processes. Precisely, in this context, the best expected cumulative reward is obtained by solving a set of equations, commonly known as Bellman's Equation. The latter is proved via Dynamic Programming Principle and, for this reason, it strongly relies on the Markov property of the environment. In the infinite horizon setting, one can show that Bellman's Equation reduces to a fixed point problem, while, in the finite horizon one, it defines a backward recursive formula. In any case, solving such an equation requires the explicit knowledge of the environment. For a more detailed explanation on decision processes we refer to \cite{bauerle2011markov}.\\
However, in Reinforcement Learning, the environment is usually not explicitly known and thus the value function has to be approximated by means of other methods. In the infinite horizon setting, the most commonly used algorithm is Watkins and Dayan's $Q$-learning algorithm (exploited in \cite{watkins1989learning,watkins1992q}). It is based on Bellman's Equation for the action-value function in its fixed point form and thus it is not directly applicable to the finite horizon setting. However, in several cases, the latter can be transformed into an infinite horizon setting, as it is shown in \cite[Section 3.4]{sutton2018reinforcement}. The convergence of the $Q$-learning algorithm has been also proven by means of Ordinary Differential Equations (ODEs) methods for recursive stochastic relations, as shown in \cite{borkar2000ode}. This approach can be also applied to obtain a method that works directly on the finite horizon setting, as done in \cite{bhatnagar2021finite}.\\
This is not the exclusive possible generalization. Indeed, the arguments described before are based on the Markov property of the environment, which is not guaranteed in general. To provide a further extension of the model, one could consider a semi-Markov decision process in place of a Markov one. These kind of processes were introduced by L\'evy in \cite{levy1954processus}, where they are defined as processes exhibiting the Markov property only on specific Markov times. In the context of decision processes, these are usually considered to be stepped (i.e. with discrete states) continuous-time semi-Markov processes, i.e. processes that can change among a certain number of discrete states via a Markov chain and such that the inter-jump times are independent of each other (while they do not need to be exponentially distributed, which is, instead, necessary and sufficient to have a stepped Markov process). One can see, for instance, \cite{khodadadi2014learning} for a model of Reinforcement Learning based on a semi-Markov decision process. See also \cite{hu2007markov} for a full description of semi-Markov decision processes.\\
In the context of Reinforcement Learning, it could be useful to consider, instead, discrete-time semi-Markov processes. These are not widely studied in literature, opposite to what happens to the continuous-time ones. semi-Markov Chains are studied, for instance, in \cite{barbu2009semi} with a matrix algebra approach. Moreover, a common strategy to construct continuous-time semi-Markov processes is through time-change (see \cite{ccinlar1972markov,ccinlar1972markovI,cinlar1974markov,jacod1974systemes} for some first considerations and, for instance, \cite{meerschaert2008triangular} for a more recent one), while this kind of approach has been considered, up to our knowledge, only in \cite{pachon2021discrete} in the discrete-time setting. Such an approach, in the continuous-time case, is revealed to be fruitful when combined with the theory of non-local operators (leading, in particular, to the generalized fractional calculus, see \cite{kochubei2011general,toaldo2015convolution}) and different properties of solutions of non-local equations are achieved via stochastic representation results (see, for instance, \cite{meerschaert2019stochastic,leonenko2013fractional,ascione2021time,ascione2021non}). This connection has not been fully explored yet in the discrete-time setting, but in \cite{pachon2021discrete} the authors achieved a first result on the link between fractional powers of the finite difference operator and the Sibuya distribution. Actually, in the continuous-time setting there is not a bijection between time-changed and semi-Markov processes, except that in the stepped case, thus the time-change approach does not cover the full range of semi-Markov processes in continuous time and space. In any case, an alternative characterization of semi-Markov process can be provided, by means of the strong Markov property of a related bivariate process involving the sojourn time. Such a characterization can be directly translated to the discrete-time setting. However, up to our knowledge, it has never been used in the construction of semi-Markov decision processes.\\
In this paper, we use the sojourn time characterization of a semi-Markov process to provide a different approach to semi-Markov decision processes. Precisely, in Section \ref{sec2} we give some preliminary notions on discrete-time semi-Markov processes. We show that all of them are obtained by a time-change procedure (thus proving that the approach presented in \cite{pachon2021discrete} can be used to describe any semi-Markov chain), by adapting a result of Jacod (in \cite{jacod1974systemes}) for stepped continuous-time semi-Markov processes. We also exploit the aforementioned characterization of a discrete-time semi-Markov process by means of sojourn time. In Section \ref{sec3} we use such a construction to provide a bivariate Markov approach to semi-Markov decision processes in which one can observe the state of the process also during the experiment. Previously, semi-Markov decision processes with continuous time allowed the agent to apply an action, that could modify the distribution of the inter-state time, only when the process changes its state (i.e. exactly where the Markov property holds). Here, we allow the agent to decide different actions while waiting for the process to change its state. Such actions can be also taken depending on the sojourn time of the environment in the current state. As a consequence, we get Bellman's Equation in terms of both the state and the sojourn time. We also prove that, in the finite horizon case, if the state space and the action space are finite, then an optimal policy exists. A numerical method for semi-Markov decision processes is exploited in Section \ref{sec4}. Such a method is based on the one shown in \cite{bhatnagar2021finite}, but it must take in consideration the fact that the set of states that can be explored depends on the present step. Finally, in Section \ref{sec5} we provide numerical experiments on two particular settings. Precisely, we first consider a problem in which only the reward explicitly depends on the sojourn time. Even in this case, the optimal policy can be chosen according to the actual sojourn time and then the resulting process is a semi-Markov one. This kind of approach could be useful to obtain an Artificial Intelligence that feels much more natural to a human observer with a \textit{gambler's fallacy} (see \cite{warren2018re}). Numerical results based on the $Q$-learning algorithm are exploited in Subsection \ref{sec512}. However, it seems to show a quadratic space complexity with respect to the horizon. To solve this problem, we also consider neural networks to provide a function approximation algorithm of $Q$: this approach is based on classical Deep Reinforcement Learning as presented in \cite{mnih2015human}. A second example is described in Subsection \ref{sec52}. This is based on a more classical problem in which semi-Markov decision processes play a prominent role. Precisely, we consider a preventive maintenance model in which the agent is partially informed, on precise time instants, of the \textit{stress level} of the system and thus can decide its \textit{regime} accordingly. Numerical results are shown in Subsection \ref{sec522}. Finally, in Section \ref{sec6} we summarize the contents of the paper and we provide some further discussion.
\section{Notions on Discrete-time semi-Markov Processes}\label{sec2}
We denote by $\N$ the set of positive integers, $\N_0=\N \cup \{0\}$ and we fix a filtered probability space $(\Omega, \cF, (\cF_t)_{t \in \N_0}, \bP)$. We also recall that, in the discrete-time and discrete-space setting, strong Markov property is implied by the Markov property (see \cite[Theorem $1.4.2$]{norris1998markov}). From now on, $E$ will denote a countable set.
\begin{definition}[{\cite[Chapter XI, Section 2]{asmussen2008applied}}]
	Let $\chi=\{\chi_n, \ n \in \N_0\}$ be a Markov chain with state space $E$. We say that $\chi$ is a \textit{jump chain} if, for any $n \in \N_0$, $\bP(\chi_{n+1}=x|\chi_{n}=x)=0$.\\
	Let $\chi$ be a jump chain and $T=\{T_n, \ n \in \N_0\}$ be a process with state space $\N$ such that $(\chi,T)=\{(\chi_n,T_n), \ n \in \N_0\}$ is a Markov process. We say that $(\chi,T)$ is a Markov additive process (MAP) if and only if
	\begin{equation*}
		\bP(\chi_{n+1}=x, T_{n+1}=t| \chi_0, \dots, \chi_n, T_0, \dots, T_n)=\bP(\chi_{n+1}=x, T_{n+1}=t| \chi_n).
	\end{equation*}
	We call $\chi$ and $T$ respectively the \textit{Markov} and \textit{additive component} of the MAP.\\
	We  say that a MAP $(\chi,T)$ is a homogeneous Markov additive processes (HMAP) $(\chi,T)$ if $\chi$ is a time-homogeneous Markov chain and for any $n,m \ge 0$
	\begin{equation*}
		\bP(\chi_{n+1}=x, T_{n+1}=t| \chi_{n})=\bP(\chi_{m+1}=x, T_{m+1}=t| \chi_m).
	\end{equation*}
\end{definition}
\begin{remark}
	One usually considers the increments of $T$ in place of $T$ when defining Markov additive processes (see, for instance, \cite[Defition $2.2$]{tang2008markov}). However, since we are working in discrete time, the formulation we provide here is equivalent.
\end{remark}
\begin{definition}[{\cite{levy1954processus}}]\label{def:semiM}
	We say that a process $X=\{X_t, \ t \in \N_0\}$ is a \textit{(discrete-time) semi-Markov process} if it exhibits the strong Markov property on any Markov time $\tau$ with values in the random set $M=\{t \in \N: \ X_{t-1} \not = X_t\}\cup \{0\}$.
\end{definition}
\begin{definition}
	Given a (discrete-time) process $X$, we define its \textit{sojourn time process} (or \textit{age process}) $\gamma=\{\gamma_t, \ t \in \N_0\}$ as follows:
	\begin{equation*}
		\gamma_t=\max\{i \le t: \ \forall j \le i, \  X_{t}=X_{t-j}\}.
	\end{equation*}
\end{definition}
Let us prove the following characterization of semi-Markov processes in terms of their sojourn time process.
\begin{theorem}\label{thm:char1}
	Let $X=\{X_t, t \in \N_0\}$ and $\gamma$ be its sojourn time process. Then the following properties are equivalent:
	\begin{enumerate}
		\item $X$ is a semi-Markov process; 
		\item $(X,\gamma)$ is a Markov process.
	\end{enumerate}
\end{theorem}
\begin{proof}
	Let us first show that $(2) \Rightarrow (1)$. Since we are working in discrete time, $(X,\gamma)$ is also a strong Markov process. Let $\tau$ be a Markov time such that $\bP(\tau \in M)=1$, where $M$ is as in Definition \ref{def:semiM}. Then we have $\bP(\gamma_\tau=0)=1$. Thus, consider $x_0,\dots,x_{\tau+1},x \in E$, with $x_0 \not = x_1$, $n \in \N$ and observe that
	\begin{align*}
		\bP(X_{\tau+n}=x|&X_\tau=x_0,X_{\tau-1}=x_1,\dots,X_{0}=x_{\tau+1})\\&=\bP(X_{\tau+n}=x|(X_\tau,\gamma_\tau)=(x_0,0),X_{\tau-1}=x_1,\dots,X_{0}=x_{\tau+1})\\
		&=\bP(X_{\tau+n}=x|(X_\tau,\gamma_\tau)=(x_0,0))=\bP(X_{\tau+n}=x|X_\tau=x_0),
	\end{align*}
	where we used the strong Markov property of $(X,\gamma)$.\\
	To prove $(1) \Rightarrow (2)$, let $n \in \N$, $x_0,\dots,x_n \in E$ and $t_0,\dots,t_n \in E$. To guarantee that we are conditioning with respect to an event with positive probability, let us suppose that for $i=1,\dots,n-1$ it holds $t_i=t_{i-1}+1$ if $x_i=x_{i-1}$ and $t_i=0$ if $x_i \not = x_{i-1}$. We want to prove that
	\begin{multline}\label{Markovprop}
		\bP((X_n,\gamma_n)=(x_n,t_n)|(X_{n-1},\gamma_{n-1})=(x_{n-1},t_{n-1}),\dots,(X_0,\gamma_0)=(x_0,t_0))\\=\bP((X_n,\gamma_n)=(x_n,t_n)|(X_{n-1},\gamma_{n-1})=(x_{n-1},t_{n-1})).
	\end{multline}
	To do this, we need to distinguish among four cases. First, let us suppose $x_n=x_{n-1}$ and $t_n \not = t_{n-1}+1$. Then, by definition of $\gamma_n$, we have
	\begin{align*}
		&\bP((X_n,\gamma_n)=(x_n,t_n)|(X_{n-1},\gamma_{n-1})=(x_{n-1},t_{n-1}),\dots,(X_0,\gamma_0)=(x_0,t_0))=0,\\
		&\bP((X_n,\gamma_n)=(x_n,t_n)|(X_{n-1},\gamma_{n-1})=(x_{n-1},t_{n-1}))=0,
	\end{align*}
	implying Equation \eqref{Markovprop}.\\
	Now let us assume $x_n=x_{n-1}$ and $t_n=t_{n-1}+1$ and rewrite
	\begin{align}\label{eq:1}
		\begin{split}
			&\bP((X_n,\gamma_n)=(x_n,t_n)|(X_{n-1},\gamma_{n-1})=(x_{n-1},t_{n-1}),\dots,(X_0,\gamma_0)=(x_0,t_0))\\&\quad =\bP(\gamma_n=t_n|X_n=x_n, (X_{n-1},\gamma_{n-1})=(x_{n-1},t_{n-1}),\dots,(X_0,\gamma_0)=(x_0,t_0))\\& \qquad \times \bP(X_n=x_n|(X_{n-1},\gamma_{n-1})=(x_{n-1},t_{n-1}),\dots,(X_0,\gamma_0)=(x_0,t_0)).
		\end{split}
	\end{align}
	Concerning the first factor of Equation \eqref{eq:1}, we have
	\begin{multline}\label{eq:2}
		\bP(\gamma_n=t_n|X_n=x_n, (X_{n-1},\gamma_{n-1})=(x_{n-1},t_{n-1}),\dots,(X_0,\gamma_0)=(x_0,t_0))=1\\=\bP(\gamma_n=t_n|X_n=x_n, (X_{n-1},\gamma_{n-1})=(x_{n-1},t_{n-1})).
	\end{multline}
	To work with the second factor of Equation \eqref{eq:1}, define the random set $M_n=\{m \le n: X_{m-1} \not = X_n\}$ and the Markov time
	\begin{equation*}
		\tau=\begin{cases} \max M_n & M_n \not = \emptyset\\
			0 & M_n=\emptyset \end{cases}
	\end{equation*}
	so that,clearly, $\bP(\tau \in M)=1$. Any event of the form $\{(X_{m},\gamma_{m})=(x_{m},t_{m})\}$ with $t-n \le m \le n-1$ is equivalent to $\{X_\tau=x_{n-1},\tau=n-t_{n-1}-1\}$, while all the other ones of the form $\{(X_{m},\gamma_{m})=(x_{m},t_{m})\}$, with $m<t-n$, are $\cF_{t-n}$-measurable. Hence, we have
	\begin{align}\label{eq:3}
		\begin{split}
			\bP&(X_n=x_n|(X_{n-1},\gamma_{n-1})=(x_{n-1},t_{n-1}),\dots,(X_0,\gamma_0)=(x_0,t_0))\\
			&=\bP(X_n=x_n|X_\tau=x_{n-1},\tau=n-t_{n-1}-1,\dots,(X_0,\gamma_0)=(x_0,t_0))\\
			&=\bP(X_{\tau+t_{n-1}+1}=x_n|X_\tau=x_{n-1},\tau=n-t_{n-1}-1,\dots,(X_0,\gamma_0)=(x_0,t_0))\\
			&=\bP(X_{\tau+t_{n-1}+1}=x_n|X_\tau=x_{n-1},\tau=n-t_{n-1}-1)\\
			&=\bP(X_{n}=x_n|(X_{n-1},\gamma_{n-1})=(x_{n-1},t_{n-1})),
		\end{split}
	\end{align}
	where we used the fact that, being $X$ semi-Markov, it exhibits the strong Markov property with respect to the jump time $\tau$. Let us stress that the previous equality holds independently of the fact that $x_n=x_{n-1}$. Combining \eqref{eq:2} and \eqref{eq:3} into \eqref{eq:1} we get
	\begin{align*}
		\begin{split}
			\bP&((X_n,\gamma_n)=(x_n,t_n)|(X_{n-1},\gamma_{n-1})=(x_{n-1},t_{n-1}),\dots,(X_0,\gamma_0)=(x_0,t_0))\\&=\bP(\gamma_n=t_n|X_n=x_n, (X_{n-1},\gamma_{n-1})=(x_{n-1},t_{n-1}),\dots,(X_0,\gamma_0)=(x_0,t_0))\\& \quad \times \bP(X_n=x_n|(X_{n-1},\gamma_{n-1})=(x_{n-1},t_{n-1}),\dots,(X_0,\gamma_0)=(x_0,t_0))\\
			&=\bP(\gamma_n=t_n|X_n=x_n, (X_{n-1},\gamma_{n-1})=(x_{n-1},t_{n-1}))\\& \quad \times \bP(X_n=x_n|(X_{n-1},\gamma_{n-1})=(x_{n-1},t_{n-1}))\\
			&=\bP((X_n,\gamma_n)=(x_n,t_n)|(X_{n-1},\gamma_{n-1})=(x_{n-1},t_{n-1})),
		\end{split}
	\end{align*}
	that is Equation \eqref{Markovprop}.\\
	Now let us consider the case $x_n \not = x_{n-1}$ and $t_n \not =0$. Then, by definition of $\gamma_n$, we have
	\begin{align*}
		&\bP((X_n,\gamma_n)=(x_n,t_n)|(X_{n-1},\gamma_{n-1})=(x_{n-1},t_{n-1}),\dots,(X_0,\gamma_0)=(x_0,t_0))=0,\\
		&\bP((X_n,\gamma_n)=(x_n,t_n)|(X_{n-1},\gamma_{n-1})=(x_{n-1},t_{n-1}))=0,
	\end{align*}
	implying Equation \eqref{Markovprop}.\\
	Finally, let us suppose that $x_n \not = x_{n-1}$ and $t_n=0$. Arguing as before, we have 
	\begin{align}\label{eq:4}
		\begin{split}
			\bP&((X_n,\gamma_n)=(x_n,0)|(X_{n-1},\gamma_{n-1})=(x_{n-1},t_{n-1}),\dots,(X_0,\gamma_0)=(x_0,t_0))\\&=\bP(\gamma_n=0|X_n=x_n, (X_{n-1},\gamma_{n-1})=(x_{n-1},t_{n-1}),\dots,(X_0,\gamma_0)=(x_0,t_0))\\& \quad \times \bP(X_n=x_n|(X_{n-1},\gamma_{n-1})=(x_{n-1},t_{n-1}),\dots,(X_0,\gamma_0)=(x_0,t_0)).
		\end{split}
	\end{align}
	Concerning the first factor of Equation \eqref{eq:4}, it holds
	\begin{multline}\label{eq:5}
		\bP(\gamma_n=0|X_n=x_n, (X_{n-1},\gamma_{n-1})=(x_{n-1},t_{n-1}),\dots,(X_0,\gamma_0)=(x_0,t_0))=1\\=\bP(\gamma_n=0|X_n=x_n, (X_{n-1},\gamma_{n-1})=(x_{n-1},t_{n-1})).
	\end{multline}
	For the second factor, let us recall, as we stated before, that Equation \eqref{eq:3} still holds. Thus, combining \eqref{eq:3} and \eqref{eq:5} into \eqref{eq:4} we obtain again \eqref{Markovprop}, concluding the proof. \qed
\end{proof}
There is a bijection between semi-Markov processes and strong Markov additive processes in the discrete-time setting. Indeed, if we start from a Markov additive process $(\chi, T)$, we can define the processes $S=\{S_n, \ n \in \N_0\}$ and $\cN=\{\cN_t, \ t \in \N_0\}$ by
\begin{align}\label{eq:N}
	\begin{split}
		S_n=\sum_{k=0}^{n}T_k, \ n \ge 0,\\
		\cN_t=\max\{n: \ S_n \le t\}
	\end{split}
\end{align}
where we set $T_0=0$. In such a case, $X_t:=\chi_{\cN_t}$ is a semi-Markov process (see \cite{ccinlar1972markovI,ccinlar1972markov,cinlar1974markov}) with sojourn time process $\gamma_t=t-S_{\cN_t}$.\\
Vice versa, if we have a semi-Markov process $X=\{X_t, \ t \in \N_0\}$, we can embed it in a continuous-time process $\widehat{X}=\{\widehat{X}_t, \ t \in [0,+\infty)\}$ by setting $\widehat{X}_t=X_{\lfloor t \rfloor}$. Let us first remark that $\widehat{X}$ is a stepped process. One can also verify that $\widehat{X}$ is semi-Markov. Indeed, if $\tau$ is a Markov time such that $\bP(\tau \in \widehat{M})=1$, where $\widehat{M}=\{t \in [0,+\infty): \ \widehat{X}_{t-}\not = \widehat{X}_t\}$, then $\tau$ is integer-valued. Moreover, $\bP(\tau \in M)=1$, where $M$ is given in Definition \ref{def:semiM}. If we consider any $n \in \N$, $x_0,\dots,x_{n+1},x \in E$, $t>0$ and $0 \le t_0<\dots<t_n<\tau$, we have
\begin{align*}
	\bP(\widehat{X}_{\tau+t}=x|&\widehat{X}_\tau=x_{n+1},\widehat{X}_{t_n}=x_n, \dots, \widehat{X}_{t_0}=x_0)\\
	&=\bP(X_{\lfloor\tau+t\rfloor}=x|X_\tau=x_{n+1},X_{\lfloor t_n \rfloor}=x_n, \dots, X_{\lfloor t_0 \rfloor}=x_0).
\end{align*}
However, being $\bP(\tau \in \N_0)=1$, it holds $\lfloor \tau+t\rfloor=\tau+\lfloor t \rfloor$ and then
\begin{align*}
	\bP(\widehat{X}_{\tau+t}=x|&\widehat{X}_\tau=x_{n+1},\widehat{X}_{t_n}=x_n, \dots, \widehat{X}_{t_0}=x_0)\\
	&=\bP(X_{\tau+\lfloor t \rfloor}=x|X_\tau=x_{n+1},X_{\lfloor t_n \rfloor}=x_n, \dots, X_{\lfloor t_0 \rfloor}=x_0)\\
	&=\bP(X_{\tau+\lfloor t \rfloor}=x|X_\tau=x_{n+1})\\
	&=\bP(\widehat{X}_{\tau+ t}=x|\widehat{X}_\tau=x_{n+1}),
\end{align*}
where we used the strong Markov property of the process $X$ on $\tau$. Thus, $\widehat{X}$ is a stepped semi-Markov process as defined by L\'evy in \cite{levy1954processus}. Hence, a result by Jacod (in \cite{jacod1974systemes}, see also \cite{cinlar1974markov} for a discussion) tells us that there exists a Markov additive process $(\chi, T)$ such that $\widehat{X}_t=\chi_{\cN_t}$ (where $\cN_t$ is defined as in \eqref{eq:N} for $t \in [0,+\infty)$). Clearly, since the jump times of $\widehat{X}$ are integer-valued random variables, $T$ has to be integer-valued and thus $\cN_t=\cN_{\lfloor t \rfloor}$. In particular, we get $X_t=\widehat{X}_{\lfloor t \rfloor}=\chi_{\cN_{\lfloor t \rfloor}}=\chi_{\cN_t}$.\\
The previous argument is the proof of the following characterization Theorem.
\begin{theorem}\label{thm:bij}
	Let $(\chi,T)$ be a Markov additive process and consider the process $\cN$ defined as in Equation \eqref{eq:N}. Then $X=\{X_t, \ t \in \N_0\}$ such that $X_t=\chi_{\cN_t}$ is a semi-Markov process. Vice versa, for any semi-Markov process $X=\{X_t, \ t \in \N_0\}$ there exists a Markov additive process $(\chi,T)$ such that $X_t=\chi_{\cN_t}$ for any $t \in \N_0$.
\end{theorem}
\begin{remark}
	In this Section we discussed some equivalent definitions of discrete-time semi-Markov processes. A further one in terms of matrix algebra arguments is given in \cite{barbu2009semi}. Moreover, let us stress that in the continuous-time case an important class of semi-Markov processes is obtained by considering $\cN$ to be independent of $\chi$. Such processes have been widely studied in literature, as they are naturally linked with generalized fractional operators. In the discrete-time case, an extensive study on the case in which the process $\cN$ is independent of $\chi$ has been carried out in \cite{pachon2021discrete}.
\end{remark}
\section{A sojourn-based approach on semi-Markov decision processes with finite horizion}\label{sec3}
Now we want to construct a decision model based on discrete-time semi-Markov processes. To do this, we will follow the line of \cite[Chapter 2]{bauerle2011markov}. Precisely, let us fix a countable \textit{state space} $E$ and the \textit{extended state space} $\widetilde{E}=E \times \N_0$, both equipped with their respective total $\sigma$-algebra. Let us consider an \textit{action space} $A$ equipped with a $\sigma$-algebra $\mathfrak{A}$ and a non-empty set-valued function $D:\widetilde{E} \to 2^A$, where, for any $(x,t) \in \widetilde{E}$, $D(x,t)$ is the set of the \textit{admissible actions}. If we equip $\widetilde{E}\times A$ with the product $\sigma$-algebra, we can define the graph $\mathcal{D}\subseteq \widetilde{E}\times A$ of $D$, i.e. $\mathcal{D}=\{(x,t,a) \in \widetilde{E} \times A: \ a \in D(x,t)\}$, and we can equip it with the trace of the aforementioned product $\sigma$-algebra. We can define a function $P$ with the following properties:
\begin{itemize}
	\item[(I)] for any $(x,t) \in \widetilde{E}$ and $a \in D(x,t)$ the function $(y,s) \in \widetilde{E} \mapsto P(y,s|x,t,a)$ is a probability mass function on $\widetilde{E}$, i.e. for any $(y,s) \in \widetilde{E}$ it holds $0 \le P(y,s|x,t,a) \le 1$ and $\sum_{(y,s) \in \widetilde{E}}P(y,s|x,t,a)=1$;
	\item[(II)] for any $(y,s) \in \widetilde{E}$ the map $(x,t,a) \in \cD \mapsto P(y,s|x,t,a)$ is measurable;
	\item[(III)] For any $x \in E$, $t,s \in \N_0$ with $s \not = t+1$ and $a \in D(x,t)$ it holds $P(x,s|x,t,a)=0$;
	\item[(IV)] For any $x,y \in E$ with $x \not = y$, $t,s \in \N_0$ with $s \not = 0$ and $a \in D(x,t)$ it holds $P(y,s|x,t,a)=0$.
\end{itemize}
We deal with the finite horizon case. Hence, let us fix our \textit{horizon} $N \in \N$. Consider \begin{equation}\label{overliner}
	(x,y,t,a) \in E \times \cD \mapsto \overline{r}_n(x|y,t,a) \in \R,\quad  n=0,\dots,N-1,
\end{equation} 
that will be our \textit{one-stage reward function}, and $(x,t) \in \widetilde{E} \mapsto g_N(x,t)$, which will be the \textit{terminal reward function}. The \textit{decision rules} will be functions $f: \widetilde{E} \to A$ such that for any $(x,t) \in \widetilde{E}$ it holds $f(x,t) \in D(x,t)$. Let us denote by $F$ the set of decision rules. A \textit{policy} $\pi$ is any element of $F^N$ (that will be denoted as $\pi=(\pi_0,\dots,\pi_{N-1})$).\\
By the Ionescu-Tulcea theorem (see, for instance, \cite{chan1974notes}), we know that for any policy $\pi$ and for any $(x,t) \in \widetilde{E}$ there exists a probability measure $\bP_{\pi,(x,t)}$ on $\Omega=\widetilde{E}^{N+1}$ equipped with the total $\sigma$-algebra (where we denote the elements as $\mathbf{x}=((x_0,t_0)\dots,(x_N,t_N))$) such that, defining $X_i(\mathbf{x})=x_i$ and $\gamma_i(\mathbf{x})=t_i$:
\begin{itemize}
	\item It holds:
	\begin{equation*}
		\bP_{\pi,(x,t)}((X_0,\gamma_0)=(y,s))=\begin{cases} 0 & (y,s)\not = (x,t) \\
			1 & (y,s)=(x,t)
		\end{cases}
	\end{equation*}
	\item For any $(y,s) \in \widetilde{E}$ and $n=0,\dots,N-1$ it holds
	\begin{equation*}
		\bP_{\pi,(x,t)}((X_{n+1},\gamma_{n+1})=(y,s)|X_0,\gamma_0,\dots,X_n,\gamma_n)=P(y,s|X_n,\gamma_n,\pi_n(X_n,\gamma_n)).
	\end{equation*}
\end{itemize}
Let us also denote by $\E_{\pi, (x,t)}$ the expected value functional associated to the measure $\bP_{\pi, (x,t)}$. As usual, the construction obtained by means of the Ionescu-Tulcea theorem guarantees that $(X,\gamma)$ is a Markov process. However, the additional properties (III) and (IV) give us some more precise structural properties of the marginal process $X$.
\begin{proposition}\label{prop:semi}
	For any policy $\pi \in F^N$ and any $x \in E$, the process $X$ is a semi-Markov process with respect to $\bP_{\pi, (x,0)}$.
\end{proposition}
\begin{proof}
	Fix $\pi \in F^N$ and $x \in E$. We actually want to show that $\gamma$ coincides, almost surely with respect to $\bP_{\pi, (x,0)}$, with the sojourn time process of $X$. Once this is done, semi-Markov property follows from the fact that $(X,\gamma)$ is a Markov process.\\
	Let us consider the sojourn time process $\widetilde{\gamma}=\{\widetilde{\gamma}_t, t \in \N_0\}$ of $X$. We want to prove that $\bP_{\pi,(x,0)}(\gamma_t=\widetilde{\gamma}_t)=1$ for any $t=0,\dots,N$. Let us first observe that for $t=0$ we have
	\begin{equation*}
		\bP_{\pi,(x,0)}(\gamma_0=\widetilde{\gamma}_0)=\bP_{\pi,(x,0)}(\gamma_0=0)=1
	\end{equation*}
	by definition of $\widetilde{\gamma}_0$ and $\bP_{\pi,(x,0)}$. Thus, let us fix $t<N$ and suppose $\bP_{\pi,(x,0)}(\gamma_t=\widetilde{\gamma}_t)=1$. This implies, in particular, that $\bP_{\pi,(x,0)}(\cdot)=\bP_{\pi,(x,0)}(\cdot|\gamma_t=\widetilde{\gamma}_t)$. We can rewrite
	\begin{equation*}
		\bP_{\pi, (x,0)}(\gamma_{t+1}=\widetilde{\gamma}_{t+1})=\E_{\pi, (x,0)}[\delta_0(\gamma_{t+1}-\widetilde{\gamma}_{t+1})],
	\end{equation*}
	where
	\begin{equation*}
		\delta_0(z)=\begin{cases} 1 & z=0 \\
			0 & z \not = 0.
		\end{cases}
	\end{equation*}
	By using the properties of the conditional expectation (see, for instance, \cite[Section $9.7$]{williams1991probability}) we get
	\begin{equation*}
		\bP_{\pi, (x,0)}(\gamma_{t+1}=\widetilde{\gamma}_{t+1})=\E_{\pi, (x,0)}[\E_{\pi, (x,0)}[\delta_0(\gamma_{t+1}-\widetilde{\gamma}_{t+1})|X_t,\gamma_t]].
	\end{equation*}
	By the Doob-Dynkin lemma (see, for instance, \cite[Lemma $2.1.24$]{bobrowski2005functional}) we know that there exists a function $f:\widetilde{E} \to \R$ such that $$\E_{\pi, (x,0)}[\delta_0(\gamma_{t+1}-\widetilde{\gamma}_{t+1})|X_t,\gamma_t]=f(X_t,\gamma_t).$$
	We can determine this function as, for any $(y,s) \in \widetilde{E}$,
	\begin{align*}
		f(y,s)&=\E_{\pi, (x,0)}[\delta_0(\gamma_{t+1}-\widetilde{\gamma}_{t+1})|X_t=y,\gamma_t=s]\\
		&=\bP_{\pi, (x,0)}(\gamma_{t+1}=\widetilde{\gamma}_{t+1}|X_t=y,\gamma_t=s)\\
		&=\bP_{\pi, (x,0)}(\gamma_{t+1}=\widetilde{\gamma}_{t+1}, X_{t+1}=y|X_t=y,\gamma_t=s)\\&+\bP_{\pi, (x,0)}(\gamma_{t+1}=\widetilde{\gamma}_{t+1}, X_{t+1} \not = y|X_t=y,\gamma_t=s)\\
		&=I_1+I_2.
	\end{align*}
	Let us first evaluate $I_1$. Observe that we are conditioning to $\gamma_t=\widetilde{\gamma}_t$ and $\gamma_t=s$, hence we also know that $\widetilde{\gamma}_t=s$. Moreover, since we know that $X_t=y$, by $\widetilde{\gamma}_t=s$ we get $X_{t-j}=y$ for any $j \le s$. With this conditioning, we have that $X_{t+1}=y$ if and only if $\widetilde{\gamma}_{t+1}=s+1$ and then
	\begin{equation*}
		I_1=\bP_{\pi, (x,0)}(\gamma_{t+1}=s+1, X_{t+1}=y|X_t=y,\gamma_t=s)=P(y,s+1|y,s,\pi_{t}(y,s)).
	\end{equation*}
	Concerning $I_2$, since we are conditioning with $X_t=y$, we know that $X_{t+1} \not = y$ if and only if $\widetilde{\gamma}_{t+1}=0$ and then
	\begin{equation*}
		I_2=\bP_{\pi, (x,0)}(\gamma_{t+1}=0, X_{t+1}\not = y|X_t=y,\gamma_t=s)=\sum_{\substack{z \in E \\ z \not = y}}P(z,0|y,s,\pi_{t}(y,s)).
	\end{equation*}
	Hence, we achieve
	\begin{align*}
		f(y,s)=P(y,s+1|y,s,\pi_{t}(y,s))+\sum_{\substack{z \in E \\ z \not = y}}P(z,0|y,s,\pi_{t}(y,s))=1,
	\end{align*}
	where we used the fact that $P$ is a probability transition function and properties (III) and (IV). The latter equality implies
	\begin{equation*}
		\bP_{\pi, (x,0)}(\gamma_{t+1}=\widetilde{\gamma}_{t+1})=\E_{\pi, (x,0)}[f(X_t,\gamma_t)]=1.
	\end{equation*}
	Finally, we conclude the proof by induction. \qed
\end{proof}
\begin{remark}
	One can actually show that, for any policy $\pi \in F^N$, $x \in E$ and $t_0 \in \N_0$, $X$ is semi-Markov with respect to $\bP_{\pi,(x,t_0)}$. However, in this case, $\gamma$ does not coincide almost surely with the sojourn time process. Precisely, if one considers the random set $A_0:=\{t \in \N_0: \ \gamma_t=0\}$,
	\begin{equation*}
		T_0:=\begin{cases} \min A_0 & A_0 \not = \emptyset \\
			+\infty & A_0=\emptyset
		\end{cases}
	\end{equation*}
	and the sojourn time process $\widetilde{\gamma}$ of $X$, then one can show that
	\begin{equation*}
		\widetilde{\gamma}_t=\begin{cases} \gamma_t & t \ge T_0 \\
			\gamma_t-t_0 & t=0,\dots,T_0-1
		\end{cases}
	\end{equation*}
	almost surely, with the same exact strategy of the previous proof. Then, one can prove that $(X,\widetilde{\gamma})$ is Markov by using the Markov property of $(X,\gamma)$, the definition of $T_0$ and properties (III) and (IV).
\end{remark}
With this idea in mind, we can refer to such model as a \textit{discrete-time semi-Markov decision process}.
\begin{remark}
	Actually, we described a Markov decision process $(X,\gamma)$ and we are asking for deterministic Markov policies. On the other hand, a semi-Markov decision process is usually described as a continuous-time decision process with non-exponential jump times, depending on the policy, which acts only on the state of the system (see, for instance, \cite{puterman2014markov}). Considering only integer-valued jump times, the bijection between Markov additive and semi-Markov processes, proved in Theorem \ref{thm:bij}, guarantees that our definition is in line with the classical one. The main difference is that we are \textit{tuning} our decisions by the observation of the sojourn times. For this reason, we refer to this approach as a \textit{sojourn-based approach}.
\end{remark}
Now let us consider the state-transition probabilities, for any $x \in E$ and any $(y,t,a) \in \cD$,
\begin{equation*}
	p(x|y,t,a)=\begin{cases} P(x,t+1|y,t,a) & x=y\\
		P(x,0|y,t,a) & x \not = y.
	\end{cases}
\end{equation*}
We can use the latter to define the \textit{one-stage expected reward function} as, for any $(y,t,a) \in \cD$,
\begin{equation*}
	r_n(y,t,a)=\sum_{x \in E}\overline{r}_n(x|y,t,a)p(x|y,t,a),
\end{equation*}
recalling that $\overline{r}_n$ is the one-stage reward function defined in \eqref{overliner}. Precisely, if $a=\pi_n(y,t)$, we get
\begin{equation*}
	r_n(y,t,\pi_n(y,t))=\E_{\pi, (x,t)}[\overline{r}_n(X_{n+1}|y,t,\pi_n(y,t)) \, | \, (X_n,\gamma_n)=(y,t)],
\end{equation*}
where $\E_{\pi, (x,t)}$ is the expected value functional associated to the measure $\bP_{\pi, (x,t)}$. In the applications, one directly defines $r_n$ in place of $\overline{r}_n$, since, as we will see in the next Sections, the function $p$ is usually unknown.\\
From now on, let us write $\bP_{\pi, n, (y,s)}=\bP_{\pi, (x,t)}(\cdot \,  | \, (X_n,\gamma_n)=(y,s))$ and $\E_{\pi, n, (y,s)}$ its expected value functional. In this case, the dependence on the initial datum is omitted for the ease of the reader, as we are going to use the aforementioned measure on events that are only involving $X_k, \gamma_k$ for $k\ge n$. In the following, we will also omit the dependence of $\bP_{\pi,(x,t)}$ from the initial datum whenever it is possible.\\
Now that we have described all the main features of the \textit{decision} process and its related rewards, let us introduce the value we want to maximize. Indeed, let us consider $\pi \in F^N$, $n \le N$ and let us define the \textit{value function} as
\begin{equation*}
	V_{n,\pi}(x,t)=\E_{\pi, n, (x,t)}\left[\sum_{k=n}^{N-1}r_k(X_k,\gamma_k, \pi_k(X_k,\gamma_k))+g_N(X_N,\gamma_N)\right]
\end{equation*}
and the \textit{optimal value function} as
\begin{equation*}
	V_n(x,t)=\sup_{\pi \in F^N}V_{n, \pi}(x,t).
\end{equation*}
%\begin{remark}
%Since we are going to work with a finite set of actions $A$, we have $V_n(x,t)=\max_{\pi \in F^N}V_{\pi,n}(x,t)$, as the set $F^N$ is finite (however, it could be very big).
%\end{remark}
We say that a policy $\pi^* \in F^N$ is \textit{optimal} if $$V_n(x,t)=V_{n, \pi^*}(x,t) \ \forall n=0,\dots, N, \ \forall (x,t) \in \widetilde{E}.$$ In practice, we are interested in finding the optimal policy $\pi^* \in F^N$ since it represents the strategy we have to adopt to maximize the expected cumulative reward $V_0$. To do this, we will rely on the dynamic programming principle. In particular, we want to obtain Bellman's Equation for $(V_0, \dots, V_N)$. To do this, let us introduce, for any $n \le N-1$, the operator $L_n$ acting on functions $v:\widetilde{E} \to \R$ as, for any $(x,t,a) \in \cD$,
\begin{equation*}
	L_n v(x,t,a)=r_n(x,t,a)+\sum_{\substack{y \in \E \\ y \not = x}}v(y,0)p(y|x,t,a)+v(x,t+1)p(x|x,t,a).
\end{equation*}
If we consider a policy $\pi \in F^N$ such that $a=\pi_n(x,t)$, then we can rewrite the previous definition as
\begin{equation*}
	L_n v(x,t,a)=\E_{\pi, n, (x,t)}[r_n(X_n,\gamma_n,\pi_n(X_n,\gamma_n))+v(X_{n+1},\gamma_{n+1})].
\end{equation*}
Moreover, for any $f \in F$ and $n \le N-1$, let us define the operator $\cT_{n,f}$ acting on functions $v:\widetilde{E} \to \R$ as, for any $(x,t) \in \widetilde{E}$,
\begin{equation*}
	\cT_{n,f}v(x,t)=L_n v(x,t,f(x,t)).
\end{equation*}
We can use the latter operator to state a Reward Iteration Theorem. It follows from \cite[Theorem $2.3.4$]{bauerle2011markov} arguing with $(X,\gamma)$ as a Markov decision process. However, for completeness, we give here the proof.
\begin{theorem}\label{thm:RI1}
	Let $\pi \in F^N$ be a policy. Then:
	\begin{itemize}
		\item[(i)] It holds $$V_{N,\pi}(x,t)=g_N(x,t)$$ for any $(x,t) \in \widetilde{E}$;
		\item[(ii)] It holds $$V_{n,\pi}(x,t)=\cT_{n,\pi_n}V_{n+1,\pi}(x,t)$$ for any $(x,t) \in \widetilde{E}$ and any $n=0,\dots,N-1$;
		\item[(iii)] It holds $$ V_{n,\pi}(x,t)=\cT_{n,\pi_n}\cdots \cT_{N-1,\pi_{N-1}}g_N(x,t)$$
		for any $(x,t) \in \widetilde{E}$ and any $n=0,\dots,N-1$.
	\end{itemize}
\end{theorem}
\begin{proof}
	Let us observe that property (i) follows from the definition of $V_{n,\pi}$ while property (iii) is clearly implied by (ii) and (i). Hence, we only have to prove property (ii). To do this, consider $n=0,\dots,N-1$, $(x,t) \in \widetilde{E}$ and evaluate $V_{n,\pi}(x,t)$. We have
	\begin{align}\label{eq:RevIt1}
		\begin{split}
			V_{n,\pi}(x,t)&=\E_{\pi, n, (x,t)}\left[\sum_{k=n}^{N-1}r_k(X_k,\gamma_k, \pi_k(X_k,\gamma_k))+g_N(X_N,\gamma_N)\right]\\
			&=r_n(x,t,\pi_n(x,t))\\
			&\qquad +\E_{\pi, n, (x,t)}\left[\sum_{k=n+1}^{N-1}r_k(X_k,\gamma_k, \pi_k(X_k,\gamma_k))+g_N(X_N,\gamma_N)\right]\\
			&=r_n(x,t,\pi_n(x,t))\\
			&\qquad +\E_{\pi, n, (x,t)}\left[\E_{\pi}\left[\sum_{k=n+1}^{N-1}r_k(X_k,\gamma_k, \pi_k(X_k,\gamma_k))+g_N(X_N,\gamma_N)\right|\left.\vphantom{\sum_{k=n+1}^{N-1}} X_{n+1},\gamma_{n+1}\right]\right],
		\end{split}
	\end{align}
	where we also used the properties of the conditional expectation. By the Doob-Dynkin lemma we know that there exists a function $f:\widetilde{E} \to \R$ such that
	$$ \E_{\pi}\left[\sum_{k=n+1}^{N-1}r_k(X_k,\gamma_k, \pi_k(X_k,\gamma_k))+g_N(X_N,\gamma_N)\right|\left.\vphantom{\sum_{k=n+1}^{N-1}} X_{n+1},\gamma_{n+1}\right]=f(X_{n+1},\gamma_{n+1}).$$
	We can evaluate $f$, for any $(x,t) \in \widetilde{E}$, as follows:
	\begin{align*}
		f(x,t)&=\E_{\pi}\left[\sum_{k=n+1}^{N-1}r_k(X_k,\gamma_k, \pi_k(X_k,\gamma_k))+g_N(X_N,\gamma_N)\right|\left.\vphantom{\sum_{k=n+1}^{N-1}} X_{n+1}=x,\gamma_{n+1}=t\right]\\
		&=\E_{\pi,n+1,(x,t)}\left[\sum_{k=n+1}^{N-1}r_k(X_k,\gamma_k, \pi_k(X_k,\gamma_k))+g_N(X_N,\gamma_N)\right]=V_{n+1, \pi}(x,t),
	\end{align*}
	where we used the Markov property of $(X,\gamma)$. Hence, equation \eqref{eq:RevIt1} becomes
	\begin{align*}
		V_{n,\pi}(x,t)&=r_n(x,t,\pi_n(x,t))+\E_{\pi, n, (x,t)}\left[V_{n+1, \pi}(X_{n+1},\gamma_{n+1})\right]\\
		&=r_n(x,t,\pi_n(x,t))+\sum_{\substack{y \in E \\ y \not = x}}V_{n+1, \pi}(y,0)p(y|x,t,\pi_n(x,t))\\
		&\qquad +V_{n+1, \pi}(x,t+1)p(x|x,t,\pi_n(x,t))\\
		&=\cT_{n,\pi_n}V_{n+1,\pi}(x,t),
	\end{align*}
	concluding the proof. \qed
\end{proof}
Now, we introduce the operator $\cT_n$, $n=0,\dots,N-1$, acting on functions $v:\widetilde{E} \to \R$, as, for any $(x,t) \in \widetilde{E}$,
\begin{equation*}
	\cT_n v(x,t)=\sup_{a \in D(x,t)} L_n v(x,t,a)=\sup_{f \in F}\cT_{n,f} v(x,t),
\end{equation*}
where the second equality follows from the fact that $F$ is the set of all selectors of $D$ and $\cT_{n,f} v(x,t)$ depends only on the value of $f$ on $(x,t)$ (and not on the whole function $f \in F$). By using this operator, we can finally state Bellman's Equation:
\begin{equation*}
	\begin{cases}
		v_N(x,t)=g_N(x,t) & (x,t) \in \widetilde{E}\\
		v_n(x,t)=\cT_n v_{n+1}(x,t) & (x,t) \in \widetilde{E}, \ n=0,\dots, N-1,
	\end{cases}
\end{equation*}
where $(v_0,\dots,v_N)$ is a $(N+1)$-uple of functions $v_i:\widetilde{E} \to \R$. Bellman's Equation clearly always admits a unique solution $(v_0,\dots,v_N)$ as each $v_n$ is defined recursively from $v_{n+1}$ and $v_N$ is set.\\
First, we want to compare the solution of Bellman's Equation with the optimal value function $(V_0,\dots,V_N)$. This is the first part of a Verification Theorem, whose proof easily follows from the fact that we reduced our decision process to a bivariate Markov decision process and from the classical Verification Theorem (see \cite[Theorem $2.3.7$]{bauerle2011markov}). Here we give the proof for completeness.
\begin{proposition}\label{prop:ver1}
	Let $(v_0,\dots,v_N)$ be the solution of Bellman's Equation. Then, for any $n=0,\dots,N$ and any $(x,t) \in \widetilde{E}$ it holds $v_n(x,t) \ge V_n(x,t)$.
\end{proposition}
\begin{proof}
	As a preliminary observation, let us notice that $L_n$ is a non-decreasing operator and thus also $\cT_{n,f}$ and $\cT_n$ share this property. \\
	Let us now observe that $V_{N,\pi}(x,t)=g_N(x,t)=v_N(x,t)$ and that, taking the supremum as $\pi \in F^N$, we get $V_N(x,t)=v_N(x,t)$. Let us now suppose that $v_{n+1}\ge V_{n+1}$. By the fact that $(v_0,\dots,v_N)$ solves Bellman's Equation and by the monotonicity of the involved operators we get
	\begin{equation*}
		v_n=\cT_n v_{n+1} \ge \cT_n V_{n+1} \ge \cT_{n, \pi_n}V_{n+1}\ge \cT_{n, \pi_n}V_{n+1,\pi}=V_{n,\pi},
	\end{equation*}
	where last equality follows from the Reward Iteration Theorem \ref{thm:RI1}. We complete the proof by (backward) induction. \qed 
\end{proof}
Now we want to state the actual Verification Theorem. To do this let us give another definition. We say that a decision rule $f \in F$ is maximal for a function $v:\widetilde{E} \to \R$ at a fixed stage $n=0,\dots,N-1$ if $\cT_n v(x,t)=\cT_{n,f} v(x,t)$ for any $(x,t) \in \widetilde{E}$. As for the previous results, the Verification Theorem easily follows from the classical one applied to the bivariate Markov process $(X,\gamma)$. We give the proof for completeness.
\begin{theorem}\label{thm:VerT}
	Let $(v_0,\dots,v_N)$ be the solution of Bellman's Equation and suppose there exists a policy $\pi^* \in F^N$ such that $\pi_n^*$ is maximal for $v_{n+1}$ at stage $n$, for any $n=0,\dots,N-1$. Then the policy $\pi^*$ is optimal and it holds, for any $n=0,\dots,N$ and $(x,t) \in \widetilde{E}$,
	\begin{equation}\label{optimal}
		V_n(x,t)=V_{n,\pi^*}(x,t)=v_n(x,t).
	\end{equation}
\end{theorem}
\begin{proof}
	Let us first observe, as before, that $v_N(x,t)=g_N(x,t)=V_N(x,t)=V_{N,\pi^*}(x,t)$ for any $(x,t) \in \widetilde{E}$. Let us suppose Equation \eqref{optimal} holds for $n+1$. We have, for any $(x,t) \in \widetilde{E}$,
	\begin{multline*}
		V_n(x,t) \le v_n(x,t)=\cT_n v_{n+1}(x,t)=\cT_{n, \pi_n^*}v_{n+1}(x,t)\\=\cT_{n, \pi_n^*}V_{n+1, \pi^*}(x,t)=V_{n, \pi^*}(x,t)\le V_n(x,t),
	\end{multline*}
	where the first inequality follows by Proposition \ref{prop:ver1} and the second last equality follows by the Reward Iteration Theorem \ref{thm:RI1}. Thus we get that Equation \eqref{optimal} holds also for $n$. We complete the proof by (backward) induction. \qed
\end{proof}
Clearly, to use the Verification Theorem to determine the optimal policy $\pi^*$, one needs to guarantee the existence of maximal functions for the solution of Bellman's Equation. In our case, it seems that even if $A$ and $E$ are finite, being $\widetilde{E}$ countable but not finite, one needs more structural assumptions to guarantee the existence of such functions. However, if we fix the horizon $N$, we can always reduce our extended state space $\widetilde{E}$ to a finite set.
\begin{proposition}\label{prop:gammabounded}
	Fix $t \in \N_0$:
	\begin{itemize}
		\item[(i)] For any $N \in \N$, $x \in E$, $\pi \in F^N$ and $n \le N$ it holds $\gamma_n \le t+n$ almost surely with respect to $\bP_{\pi, (x,t)}$;
		\item[(ii)] For any $N \in \N$, $x \in E$, $\pi \in F^N$ and $n=0,\dots,N$  it holds $\gamma_n \le t+N$ almost surely with respect to $\bP_{\pi, (x,t)}$.
	\end{itemize}
\end{proposition}
\begin{proof}
	Clearly, property (ii) follows from (i), so we only need to prove (i). To do this, fix any $x \in E$ and $\pi \in F^N$. Let us first observe that $\bP_{\pi, (x,t)}(\gamma_0=t)=1$. Let us now suppose that $\bP_{\pi, (x,t)}(\gamma_n \le t+n)=1$ for some $n<N$. Consider any $j \ge 2$ and evaluate $\bP_{\pi, (x,t)}(\gamma_{n+1}=t+n+j)$. To do this, let us rewrite
	\begin{equation*}
		\bP_{\pi, (x,t)}(\gamma_{n+1}=t+n+j)=\E_{\pi, (x,t)}[\delta_{0}(\gamma_{n+1}-(t+n+j))].
	\end{equation*}
	By the properties of conditional expectation we have
	\begin{equation*}
		\bP_{\pi, (x,t)}(\gamma_{n+1}=t+n+j)=\E_{\pi, (x,t)}[\E_{\pi, (x,t)}[\delta_{0}(\gamma_{n+1}-(t+n+j))|X_n, \gamma_n]].
	\end{equation*}
	By the Doob-Dynkin lemma we can set, for any $(y,s) \in \widetilde{E}$, $$f_j(y,s)=\E_{\pi, (x,t)}[\delta_{0}(\gamma_{n+1}-(t+n+j))|X_n=y, \gamma_n=s]$$ to get
	\begin{equation*}
		\bP_{\pi, (x,t)}(\gamma_{n+1}=t+n+j)=\E_{\pi, (x,t)}[f_j(X_n,\gamma_n)]=\E_{\pi, (x,t)}[f_j(X_n,\gamma_n)\mathbf{1}_{\{0,\dots,t+n\}}(\gamma_n)],
	\end{equation*}
	where
	$$\mathbf{1}_{\{0,\dots,t+n\}}(z)=\begin{cases} 1 & z \in \{0,\dots,t+n\} \\ 0 & z \not \in \{0,\dots,t+n\},\end{cases}$$
	and we used the fact that $\bP_{\pi, (x,t)}(\gamma_n \le t+n)=1$. Thus, let us fix $y \in E$, $s \le t+n$ and let us evaluate $f_j(y,s)$. We have
	\begin{align*}
		f_j(y,s)&=\E_{\pi, (x,t)}[\delta_{0}(\gamma_{n+1}-(t+n+j))|X_n=y, \gamma_n=s]\\
		&=\sum_{z \in E}\bP_{\pi, (x,t)}(\gamma_{n+1}=t+n+j,X_{n+1}=z|X_n=y, \gamma_n=s)\\
		&=\sum_{\substack{z \in E \\ z \not = y}}P(z,t+n+j|y,s,\pi_n(y,s))+P(y,t+n+j|y,s,\pi_n(y,s)).
	\end{align*}
	Now observe that $P(z,t+n+j|y,s,\pi_n(y,s))=0$ for any $z \not = y$ by property (IV) and $t+n+j \ge 2>0$. On the other hand, we also have $t+n+j\ge t+n+2 >t+n+1 \ge s+1$ and then, by property (III), we get $P(y,t+n+j|y,s,\pi_n(y,s))=0$. Thus we conclude that
	\begin{equation*}
		f_j(y,s)=0, \ \forall y \in E, \ \forall s \le t+n
	\end{equation*}
	and
	\begin{equation*}
		\bP_{\pi, (x,t)}(\gamma_{n+1}=t+n+j)=0.
	\end{equation*}
	Finally we achieve
	\begin{multline*}
		\bP_{\pi, (x,t)}(\gamma_{n+1}\le t+n+1)= 1-\bP_{\pi, (x,t)}(\gamma_{n+1}> t+n+1)\\=1-\sum_{j=2}^{+\infty}\bP_{\pi, (x,t)}(\gamma_{n+1}=t+n+j)=1.
	\end{multline*}
	We complete the proof by induction. \qed 
\end{proof}
The latter Proposition guarantees that, supposed we fixed the horizon $N$ and the initial state $(x,t)$, we can substitute the extended state space $\widetilde{E}$ with $\widetilde{E}_{N,t}=E \times \{0,\dots,N+t\}$. Since we will work in general with the initial state $(x,0)$ (so that $\gamma$ is actually the sojourn time process), let us denote $\widetilde{E}_N=E \times \{0,\dots,N\}$.
We can use this observation to prove an Existence Theorem, at least in the simplest case.
\begin{theorem}\label{thm:existence}
	Let $E$ and $A$ be finite sets (equipped with the total $\sigma$-algebra). Fix any horizon $N$ and suppose that for any $(y,s) \not \in \widetilde{E}_N$ it holds $g_N(y,s)=0$ and $r_n(y,s,a)=0$ for $n=0,\dots,N-1$. Then there exists an optimal policy $\pi^* \in F^N$.
\end{theorem}
\begin{proof}
	Let $(v_0,\dots,v_N)$ be the solution of Bellman's Equation. We only need to show that for any $n=0,\dots,N-1$ there exists a decision rule $\pi_n^*$ that is maximal for $v_{n+1}$ at stage $n$. First of all, let us observe that $v_N(y,s)=g_N(y,s)=0$ whenever $(y,s) \not \in \widetilde{E}_N$. Next, let us consider $n=N-1$. Recall that
	\begin{equation*}
		v_{N-1}(y,s)=\cT_{N-1}v_{N}(y,s)=\sup_{f \in F}\cT_{N-1,f}v_{N}(y,s).
	\end{equation*}
	Fix any decision rule $f \in F$ and observe that, for $(y,s) \not \in \widetilde{E}_N$, it holds
	\begin{equation*}
		\cT_{N-1,f}v_{N}(y,s)=r_N(y,s,f(y,s))+v_{N}(y,s)=0,
	\end{equation*}
	hence also $v_{N-1}(y,s)=0$ for any $(y,s) \not \in \widetilde{E}_N$. Let us now consider two functions $f_1,f_2 \in F$ such that for any $(y,s) \in \widetilde{E}_N$ it holds $f_1(y,s)=f_2(y,s)$. Then we have, for $(y,s) \in \widetilde{E}_N$,
	\begin{multline*}
		\cT_{N-1,f_1}v_{N}(y,s)=r_n(y,s,f_1(y,s))+v_N(y,s)\\=r_n(y,s,f_2(y,s))+v_N(y,s)=\cT_{N-1,f_2}v_{N}(y,s),
	\end{multline*}
	while, for $(y,s) \not \in \widetilde{E}_N$,
	\begin{equation*}
		\cT_{N-1,f_1}v_{N}(y,s)=0=\fT_{N-1,f_2}v_{N}(y,s),
	\end{equation*}
	as observed before.
	Let us introduce the function $\Phi_{N-1}: f \in F \mapsto \cT_{N-1,f}v_{N}\in \R^{\widetilde{E}}$ (where $\R^{\widetilde{E}}$ is the set of functions that map $\widetilde{E}$ to $\R$). We can also define on $F$ the equivalence relation $\sim$ as follows: for two decision rules $f_1,f_2 \in F$ we say that $f_1 \sim f_2$ if and only if $f_1(y,s)=f_2(y,s)$ for any $(y,s) \in \widetilde{E}_N$. Let us denote by $[f]$ the equivalence class of $f$, $F/\sim$ the quotient set and $\mathfrak{p}:F \to F/\sim$ the projection. Since we have $\Phi_{N-1}(f_1)=\Phi_{N-1}(f_2)$ whenever $f_1 \sim f_2$, there exists a unique function $\widetilde{\Phi}_{N-1}:F/\sim \to \R^{\widetilde{E}}$ such that for any $f \in F$ it holds $\widetilde{\Phi}_{N-1}(\mathfrak{p}(f))=\Phi_{N-1}(f)$. Hence we can rewrite
	\begin{equation*}
		\cT_{N-1}v_{N}=\sup_{[f] \in F/\sim}\widetilde{\Phi}_{N-1}([f]).
	\end{equation*}
	However, it is not difficult to check that there is a bijection between $F/\sim$ and the set of functions mapping $\widetilde{E}_N \to A$, that is finite since both $\widetilde{E}_N$ and $A$ are finite. Hence there exists $[f_{N-1}^*]$ such that
	\begin{equation*}
		\cT_{N-1}v_{N}=\widetilde{\Phi}_{N-1}([f_{N-1}^*]).
	\end{equation*}
	Any decision rule $\pi_{N-1}^* \in [f_{N-1}^*]$ is maximal for $v_{N}$ at stage $N-1$.\\
	Now let us suppose that we have shown the existence of a maximal decision rule $\pi_{n}^*$ for $v_{n+1}$ at stage $n$ for some $n=1,\dots,N-1$ and let us suppose that $v_n(y,s)=0$ for any $(y,s) \not \in \widetilde{E}$. First, observing that for any $f \in F$ and $(y,s) \not \in \widetilde{E}$ it holds
	\begin{equation*}
		\cT_{n-1,f}v_{n}(y,s)=r_{n-1}(y,s,f(y,s))+v_{n}(y,s)=0,
	\end{equation*}
	we get $v_{n-1}(y,s)=0$ for any $(y,s) \not \in \widetilde{E}_N$. Moreover, with the same exact argument as before, if we define the map $\Phi_{n-1}:f\in F \to \cT_{n-1,f}v_{n} \in \R^{\widetilde{E}}$, it holds $\Phi_{n-1}(f_1)=\Phi_{n-1}(f_2)$ whenever $f_1 \sim f_2$. Hence there exists a unique function $\widetilde{\Phi}_{n-1}:F/\sim \to \R^{\widetilde{E}}$ such that, for any $f \in F$, $\widetilde{\Phi}_{n-1}(\mathfrak{p}(f))=\Phi_{n-1}(f)$. Again, we have
	\begin{equation*}
		\cT_{n-1}v_{n}=\sup_{[f] \in F/\sim}\widetilde{\Phi}_{n-1}([f])
	\end{equation*}
	thus, since $F/\sim$ is finite, we can find $[f_{n-1}^*]$ such that
	\begin{equation*}
		\cT_{n-1}v_{n}=\widetilde{\Phi}_{n-1}([f_{n-1}^*]).
	\end{equation*}
	Any decision rule $\pi_{n-1}^* \in [f_{n-1}^*]$ is maximal for $v_{n}$ at stage $n-1$. We complete the proof by induction. \qed
\end{proof}
In the following, we need another auxiliary function that is strictly linked with the value function. Precisely, we define the \textit{quality function} or \textit{action-value function} (see \cite{sutton2018reinforcement}) as
\begin{align*}
	Q_{n,\pi}:(x,t,a) \in \cD \mapsto \E_{\pi,n,(x,t)}&\left[\sum_{k=n+1}^{N-1}r_k(X_k,\gamma_k,\pi(X_k,\gamma_k))\right.\\
	&\left.\vphantom{\sum_{k=n+1}^{N-1}}\qquad +g_N(X_N,\gamma_N)+r_n(X_n,\gamma_n,a)\right].
\end{align*}
One can also define the \textit{optimal action-value function} as
\begin{equation*}
	Q_{n}^*(x,t,a)=\sup_{\pi \in F^N}Q_{n,\pi}(x,t,a), \ \forall (x,t,a) \in \cD
\end{equation*}
so that, by definition,
\begin{equation*}
	V_n(x,t)=\sup_{a \in A}Q_{n}^*(x,t,a), \ \forall (x,t) \in \widetilde{E}.
\end{equation*}
Clearly, the previous arguments apply also on $Q$. Precisely, let us define the operators $\widehat{\cT}_{n,f}$, for $n=0,\dots,N-1$ and $f \in F$, acting on functions $v:\cD \to \R$ as, for any $(x,t,a) \in \cD$,
\begin{multline*}
	\widehat{\cT}_{n,f}v(x,t,a)=r_n(x,t,a)+\sum_{\substack{y \in E \\ y \not = x}}v(y,0,f(y,0))p(y|x,t,a)+v(x,t+1,f(x,t+1))p(x|x,t,a)
\end{multline*}
and $\widehat{\cT}_n$, for $n=0,\dots,N-1$, acting on functions $v:\cD \to \R$ as, for any $(x,t,a) \in \cD$,
\begin{equation*}
	\widehat{\cT}_n v(x,t,a)=\sup_{f \in F}\widehat{\cT}_{n,f}v(x,t,a).
\end{equation*}
Again, we say that $f \in F$ is maximal for a function $v:\cD \to \R$ for a fixed stage $n=0,\dots,N-1$ if $\widehat{\cT}_n v=\widehat{\cT}_{n,f}v$. With this in mind, we get the following result.
\begin{theorem}\label{thm:theoremQ}
	The following properties are true:
	\begin{itemize}
		\item[(i)] For any $\pi \in F^N$ it holds
		\begin{equation*}
			Q_{N,\pi}(x,t,a)=g_N(x,t), \qquad \forall (x,t,a) \in \cD
		\end{equation*}
		\item[(ii)] For any $\pi \in F^N$ it holds
		\begin{equation*}
			Q_{n,\pi}(x,t,a)=\widehat{\cT}_{n,\pi_{n+1}}Q_{n+1,\pi}(x,t,a), \qquad \forall (x,t,a) \in \cD
		\end{equation*}
		for any $n=0,\dots,N-1$
		\item[(iii)] Let $(q_0,\dots,q_N)$ be the solution of
		\begin{equation}\label{eq:BEqn}
			\begin{cases}
				q_N(x,t,a)=g_N(x,t) & \forall (x,t,a) \in \cD\\
				q_n(x,t,a)=\widehat{\cT}_n q_{n+1}(x,t,a) & \forall (x,t,a) \in \cD, \ n=0,\dots,N-1.
			\end{cases}
		\end{equation}
		Then it holds
		\begin{equation*}
			q_n(x,t,a) \ge Q_n^*(x,t,a) \qquad \forall (x,t,a) \in \widetilde{E} \times A.
		\end{equation*}
		\item[(iv)] Let $(q_0,\dots,q_N)$ be the solution of Equation \eqref{eq:BEqn} and suppose there exists a policy $\pi^* \in F^N$ such that $\pi_{n+1}^*$ is maximal for $q_n$ at stage $n$ for any $n=0,\dots,N-2$. Then $Q_{n,\pi^*}=Q_n^*=q_n$ for any $n=0,\dots,N$.
	\end{itemize}
\end{theorem}
We omit the proof of the previous theorem as it is identical to the ones we presented before for the value function $V_n$. Moreover, it is clear that the policy $\pi^*$ is optimal and that $V_n(x,t)=Q_{n,\pi^*}(x,t,\pi^*_n(x,t))$. We will work with the action-value function $Q$ in what follows.\\
Before proceeding to the next section, let us summarize what we did:
\begin{itemize}
	\item We set the notation and introduced the decision process, thanks to the Ionescu-Tulcea theorem, and we proved, in Proposition \ref{prop:semi}, that such a process is semi-Markov;
	\item Once this is done, we introduced the one-stage expected reward functions, the value functions and the operators that will lead to Bellman's Equation;
	\item A Reward Iteration Theorem (Theorem \ref{thm:RI1}) is proved, to show a backward recursive method from which one deduces the value function from the policy (via the operators $\cT_{n,f}$) and the reward functions (both the one-stage ones and the terminal one);
	\item Another operator is then introduced, by taking the supremum on $\cT_{n,f}$ in the variable $f \in F$, which is used to define a backward recursive formula, which is Bellman's Equation;
	\item Thanks to the Reward Iteration Theorem, we prove a Verification Theorem, whose proof is split in Proposition \ref{prop:ver1} and Theorem \ref{thm:VerT};
	\item The Verification Theorem tells us how to determine the optimal policy if it exists: one has to solve Bellman's Equation and then try to determine the maximal decision $\pi^*_n$, such that $v_n=\cT_n v_{n+1}=\cT_{n, \pi^*_n} v_{n+1}$;
	\item Finally, we need to prove that such an optimal decision exists, so that the previous methodology can be applied: a criterion for existence is proved in Theorem \ref{thm:existence}, which also guarantees that one has to search for $\pi^*_n$ in a finite set of decisions;
	\item The same results are then extended to the quality function in Theorem \ref{thm:theoremQ}.
\end{itemize}
Thus, \textit{in theory}, we could find an exact solution of the decision problem. However, we usually lack information on the environment. Such a problem is investigated in the next section.
\section{Sojourn-based semi-Markov Reinforcement Learning}\label{sec4}
Let us stress out that, if we want to work with a discrete-time semi-Markov decision process by directly using Bellman's Equation, we need to know not only the (parametric) transition probability matrix of the jump chain, but also the (parametric) distribution of the inter-jump times. In our approach, this coincides with the estimate of the parametric transition probability (multidimensional) matrix $P(y,s|x,t,a)$ of the bivariate process $(X,\gamma)$, recalling that it is subject to both properties (III) and (IV). In any case, such information is usually unknown to the agent. Thus, estimation has to be carried on with a different strategy.\\
For Markov decision processes with infinite horizon, the most famous algorithm to estimate $Q^*$ is the $Q$-learning algorithm (whose convergence was first proved in \cite{watkins1992q}). However, such an algorithm is based, in particular, on the fact that, for stationary rewards, $Q^*$ (or, equivalently, $V$) does not depend on the step $n$ (see, for instance, the Verification Theorem given in \cite[Theorem 7.1.7]{bauerle2011markov}) and that Bellman's Equation reduces to a fixed point relation. This is not the case in the finite horizon setting.\\
Recently, in \cite{bhatnagar2021finite}, the authors provided an algorithm that deals with the finite-horizon case. Let us adapt their results to our model. For a fixed horizon $N$, we define (possibly randomly) the \textit{initial $Q$-iterate} $Q^0_n(x,t,a)$ for $n=0,\dots,N-1$ and we set $Q^0_N(x,t,a)=g_N(x,t)$. Then, supposed $Q^m_n(x,t,a)$ is defined for any $n=1,\dots,N$, we set $Q^{m+1}_N(x,t,a)=g_N(x,y)$ and
\begin{multline}\label{Qlearn}
	Q^{m+1}_n(x,t,a)=(1-\alpha(m))Q^m_n(x,t,a)+\alpha(m)(r_n(x,t,a)\\+\max_{b \in D(X_{n+1},\gamma_{n+1})}Q^{m}_{n+1}(X^{m}_{n+1},\gamma^{m}_{n+1},b)), \quad n=0,\dots,N-1,
\end{multline}
where $\alpha:\N_0 \to [0,1]$ and $(X^{m}_{n+1},\gamma^{m}_{n+1})$ is obtained by sampling when $X^{m}_n=x$, $\gamma^{m}_n=t$ and the action $a$ is taken at step $n$. We also assume that, for any $m \not = m^\prime$, $(X^{m}_{n},\gamma^m_{n})_{n=0,\dots,N}$ is independent of $(X^{m^\prime}_{n},\gamma^{m^\prime}_{n})_{n=0,\dots,N}$. Moreover, we ask that, for any $m \in \N$, $(X^{m}_{0},\gamma^m_{0})$ is distributed according to a fixed probability law $P_0: \widetilde{E}_N \to [0,1]$ such that for any $t \not = 0$ it holds $P_0(x,t)=0$. Let us stress that such an update is carried on only for triples $(x,t,a)$ such that $(x,t)=(X^m_n,\gamma^m_n)$ and $a$ is selected in $\arg \max Q^{m}_{n}(X^{m}_{n},\gamma^{m}_{n},\cdot)$, while for all the other triples one sets $Q^{m+1}_n(x,t,a)=Q^{m}_n(x,t,a)$. \\
The quantity $\alpha(m)$ is called the \textit{learning rate} of the recursive relation. Clearly, if we start from a point of the form $(x,0)$, we cannot visit all the states of the form $(y,t)$ for any $t=0,\dots,N$, but only the states that belong to $\widetilde{E}_N$. Thus, we have to study the functions $Q^m_n$ only on the set $\widetilde{E}_N$. Moreover, if we are considering the step $n$, Proposition \ref{prop:gammabounded} tells us that we can only consider the states of the form $(x,t)$ with $t \le n$. Thus, let us consider $\widetilde{E}_n$ as the \textit{reachable states spaces at step $n$}, recalling that
\begin{equation*}
	\widetilde{E}_{n}=\{(x,t) \in E \times \N_0: \ t \le n\}.
\end{equation*}
The idea here is to incorporate the step $n$ as an argument of the quality function. To do this, we define the \textit{reachable step-state space} of the tri-variate Markov process $(n,X_n,\gamma_n)$ as
$$\overline{E}_N:=\{(n,x,t) \in \{0,\dots,N\} \times \widetilde{E}_N: \ (x,t) \in \widetilde{E}_n\}$$
and the \textit{admissible step-state-action space}
$$\overline{\cD}_N:=\{(n,x,t,a) \in \{0,\dots,N\} \times \widetilde{E}_N: \ (x,t) \in \widetilde{E}_n, \  a \in D(x,t)\}.$$
Once this is done, we can define the functions $\mathbf{Q}^m:\overline{\cD}_N \mapsto \R$ as
\begin{equation*}
	\mathbf{Q}^m(n,x,t,a)=Q_n^m(x,t,a), \quad \forall (n,x,t,a) \in \overline{\cD}_N, \ m \in \N_0 
\end{equation*}
and $\mathbf{Q}^*:\overline{\cD}_N \mapsto \R$ as
\begin{equation*}
	\mathbf{Q}^*(n,x,t,a)=Q_n^*(x,t,a) \quad  \forall (n,x,t,a) \in \overline{\cD}_N.
\end{equation*}
With these notations, we can rewrite equation \eqref{Qlearn} in terms of the functions $\mathbf{Q}^m$:
\begin{multline}\label{Qlearnfunctional}
	\mathbf{Q}^{m+1}(n,x,t,a)=\mathbf{Q}^m(n,x,t,a)+\alpha(m)(\mathbf{r}(n,x,t,a)\\+\max_{b \in D(X^{m+1}_{n+1},\gamma^{m+1}_{n+1})}\mathbf{Q}^m(n+1,X_{n+1}^{m+1},\gamma_{n+1}^{m+1},b)-\mathbf{Q}^m(n,x,t,a)), \\ (n,x,t,a) \in \overline{\cD}_N \mbox{ with } n=0,\dots,N-1
\end{multline}
where $\mathbf{r}:\overline{\cD}_N \to \R$ is defined as
\begin{equation*}
	\mathbf{r}(n,x,t,a)=\begin{cases} r_n(n,x,t,a) & n=0,\dots,N-1 \\
		g_N(x,t) & n=N. \end{cases}
\end{equation*}
Moreover, we set $\mathbf{Q}^m(N,x,t,a)=g_N(x,t)$ for any $(N,x,t,a) \in \overline{\cD}_N$. Since $\overline{\cD}_N$ is finite, we can identify any function $\mathbf{Q}: \overline{\cD}_N \to \R$ with a vector in $\R^{|\overline{\cD}_N|}$, where $|\overline{\cD}_N|$ is the cardinality of $\overline{\cD}_N$. However, to avoid ambiguity of notations, due to possible different indexing of vectors in $\R^{|\overline{\cD}_N|}$, we denote the space of such functions as $\overline{F}_N$, equipped with the norm
\begin{equation*}
	\Norm{\mathbf{Q}}{\infty}:=\max_{(n,x,t,a) \in \overline{\cD}_N}|\mathbf{Q}(n,x,t,a)|, \quad \forall \mathbf{Q} \in \overline{F}_N.
\end{equation*}
Formula \eqref{Qlearnfunctional} tells us that $\mathbf{Q}^{m+1}$ is constructed from $\mathbf{Q}^m$ by resorting to stochastic sampling. For such a reason, $\mathbf{Q}^m$ should be a random variable on a suitable probability space. To construct it, let us assume that $\mathbf{Q}^0$ is initialized randomly according to the probability law $\overline{P}_0$ on $(\overline{F}_N, \cB(\overline{F}_N))$, where $\cB(\overline{F}_N)$ is the Borel $\sigma$-algebra of $\overline{F}_N$. Next, we consider the set
\begin{equation*}
	\widetilde{E}_N^{(N)}=\{(x_n,t_n)_{n \le N} \in (\widetilde{E}_N)^N: \ (x_n,t_n) \in \widetilde{E}_n,\ \forall n \le N\}
\end{equation*}
and, for any $\mathbf{Q} \in \overline{F}_N$ and $m \in \N$, we define the function $F_{\mathbf{Q},m}:\widetilde{E}_N^{(N)} \to \overline{F}_N$ by asking that, for any $(\mathbf{x},\mathbf{t}) \in \widetilde{E}_N^{(N)}$, it holds
\begin{align*}
	(F_{\mathbf{Q},m}(\mathbf{x},\mathbf{t}))&(n,y,s,a)\\&=\mathbf{Q}(n,y,s,a)+\alpha(m)(\mathbf{r}(n,y,s,a)\\
	&+\max_{b \in D(x_{n+1},t_{n+1})}\mathbf{Q}(n+1,x_{n+1},t_{n+1},b)-\mathbf{Q}(n,y,s,a))
\end{align*}
for any $(n,y,s,a) \in \overline{\cD}_N$ such that $n \le N-1$, while $(F_{\mathbf{Q}}(\mathbf{x},\mathbf{t}))(N,y,s,a)=g_N(y,s)$ for any $(N,y,s,a) \in \overline{\cD}_N$. Let $\cB(\overline{F}_N)$ be the Borel $\sigma$-algebra of $\overline{F}_N$. For any $m \in \N_0$ we define the function $\overline{P}_m:\cB(\overline{F}_N) \times \overline{F}_N \mapsto [0,1]$ as:
\begin{equation*}
	\overline{P}_{m+1}(B|\mathbf{Q})=\sum_{x \in E}\bP_{\pi^{m},(x,0)}(F_{\mathbf{Q},m}(X^{m},\gamma^{m}) \in B)P_0(x,0),
\end{equation*}
where $\pi^m \in F^N$ is any policy such that, for any $(n,x,t) \in\overline{E}_N$, it holds $\pi^m_n(x,t) \in \arg \max \mathbf{Q}(n,x,t)$. For fixed $\mathbf{Q} \in \overline{F}_N$, $\overline{P}_{m+1}(\cdot |\mathbf{Q})$ is a probability measure on $(\overline{F}_N, \cB(\overline{F}_N))$, while, for fixed $B \in \cB(\overline{F}_N)$, $\overline{P}_{m+1}(B |\cdot)$ is a measurable function. Thus, by the Ionescu-Tulcea theorem, we can construct a probability space $(\overline{\Omega}, \overline{\cF}, \bP)$ supporting the sequence $(\mathbf{Q}^m)_{m \in \N}$ defined via the recursive relation \eqref{Qlearnfunctional} such that:
\begin{itemize}
	\item $\overline{\Omega}=\prod_{m=0}^{\infty}\overline{F}_N$, $\overline{\cF}=\otimes_{m=0}^{\infty}\cB(\overline{F}_N)$ and the process $(\mathbf{Q}^m)_{m \in \N}$ is Markov;
	\item For any $B \in \cB(\overline{F}_N)$ it holds
	\begin{equation*}
		\bP(\mathbf{Q}^{m+1} \in B|\mathbf{Q}^m)=\overline{P}_{m+1}(B|\mathbf{Q});
	\end{equation*}
	\item For any $B \in \cB(\overline{F}_N)$ it holds
	\begin{equation*}
		\bP(\mathbf{Q}^{0} \in B)=\overline{P}_0(B).
	\end{equation*}
\end{itemize}
Once this construction is done, we can study the convergence of such a stochastic approximation algorithm. To do this, let us first introduce the \textit{absolute error}
\begin{equation*}
	\mathbf{e}(m)=\Norm{\mathbf{Q}^m-\mathbf{Q}^*}{\infty} \quad m \in \N_0.
\end{equation*}
Now, to prove any convergence result, we will make use (as in \cite{bhatnagar2021finite}) of the ODE method presented in \cite{borkar2000ode} for stochastic approximation algorithms. However, we first need to recast formula \eqref{Qlearnfunctional}. Precisely, we add and subtract $$\alpha(m)\E_{\pi^{(a,n)},n,(x,t)}\left[\max_{b \in D(X^{m}_{n+1},\gamma^{m}_{n+1})}\mathbf{Q}^m(n+1,X_{n+1}^{m},\gamma_{n+1}^{m},b)\right],$$
where with $\pi^{(a,n)}$ we denote any policy $\pi$ such that $\pi_{n}(x,t)=a$. Let us also set $\mathbf{M}^{m+1}:\overline{\cD}_N \to \R$ as 
\begin{align}\label{eq:M}
	\begin{split}
		\mathbf{M}^{m+1}(n,x,t,a)&=\max_{b \in D(X^{m}_{n+1},\gamma^{m}_{n+1})}\mathbf{Q}^m(n+1,X_{n+1}^{m},\gamma_{n+1}^{m},b)\\
		&-\E_{\pi^{(a,n)},n,(x,t)}\left[\left.\max_{b \in D(X^{m}_{n+1},\gamma^{m}_{n+1})}\mathbf{Q}^m(n+1,X_{n+1}^{m},\gamma_{n+1}^{m},b)\right|\mathbf{Q}^m \right]
	\end{split}
\end{align}
if $n \not = N$, while $\mathbf{M}^{m+1}(N,x,t,a)=0$,
and define the operator $T:\overline{F}_N \to \overline{F}_N$ acting on $\mathbf{Q}:\overline{\cD}_N \to \R$ as
\begin{multline*}
	T(\mathbf{Q})(n,x,t,a)=\mathbf{r}(n,x,t,a)\\+\E_{\pi^{(a,n)},n,(x,t)}\left[\max_{b \in D(X^{m}_{n+1},\gamma^{m}_{n+1})}\mathbf{Q}(n+1,X_{n+1}^{m},\gamma_{n+1}^{m},b)\right]
\end{multline*}
if $n \not = N$ and $T(\mathbf{Q})(N,x,t,a)=\mathbf{Q}(N,x,t,a)$. Let us stress that the definition of $T$ is independent of $m$ thanks to the fact that the sequences $(X^{m}_n,\gamma^{m}_n)_{n \le N}$ are independent of each other and identically distributed once the policy is prescribed. Finally, if we also introduce $h(\mathbf{Q})=T(\mathbf{Q})-\mathbf{Q}$, we can rewrite the relation \eqref{Qlearnfunctional} as
\begin{multline}\label{Qlearnfunctional2}
	\mathbf{Q}^{m+1}(n,x,t,a)=\mathbf{Q}^m(n,x,t,a)+\alpha(m)(h(\mathbf{Q}^{m})(n,x,t,a)+\mathbf{M}^{m+1}(n,x,t,a)),
\end{multline}
for any $(n,x,t,a) \in \overline{\cD}_N$. Let us stress that, since $\overline{\cD}_N$ is a finite set, $\mathbf{Q}^m$, $\mathbf{Q}^*$ and $\mathbf{M}^m$ can be considered as random variables with values in $\R^{|\overline{\cD}_N|}$. To fix the notation, let us suppose that we represent them as column vectors. By definition of $\mathbf{M}^m$, one easily gets $\E[\mathbf{M}^m]=0$, while, since $\{(X^{m}_n,\gamma^{m}_n)_{n=0,\dots,N}\}_{m \in \N}$ are independent of each other, for any $m \not = m^\prime$ it holds $\E[\mathbf{M}^m(\mathbf{M}^{m^\prime})^T]=0$. Thus, Equation \eqref{Qlearnfunctional2} can be seen as a recursive stochastic approximation algorithm of the form discussed in \cite{borkar2000ode}. Up to considering the bivariate process $(X_n,\gamma_n)$, this is the strategy adopted in \cite{bhatnagar2021finite} to prove the convergence of the $Q$-learning algorithm. Let us give here some details for completeness. To do this, we need the following Lemma, which is shown inside the proof of \cite[Proposition 2]{bhatnagar2021finite}. We give here an alternative proof in the case of finite sets.
\begin{lemma}\label{lemma:max}
	Let $\cS$ be any finite set and consider two functions $F_i: \overline{\cD}_N \to \R$, $i=1,2$. Then it holds
	\begin{equation}\label{eq:maxrel}
		\left|\max_{x \in \cS}F_1(x)-\max_{x \in \cS}F_2(x)\right| \le \max_{x \in \cS}|F_1(x)-F_2(x)|.
	\end{equation}
\end{lemma}
\begin{proof}
	Let $x_i \in \overline{\cD}_N$, $i=1,2$, be such that
	\begin{equation*}
		F_i(x_i)=\max_{x \in \cS}F_i(x).
	\end{equation*}
	If $F_1(x_1)=F_2(x_2)$, then inequality \eqref{eq:maxrel} clearly holds. Thus let us suppose, without loss of generality, that $F_1(x_1)>F_2(x_2)$. By definition of $x_2$ we have $F_2(x_2) \ge F_2(x_1)$ and then
	\begin{align*}
		\left|\max_{x \in \cS}F_1(x)-\max_{x \in \cS}F_2(x)\right|&= F_1(x_1)-F_2(x_2)\\ &\le  F_1(x_1)-F_2(x_1)\le \max_{x \in \cS}|F_1(x)-F_2(x)|,
	\end{align*}
	concluding the proof. \qed
\end{proof}
As a consequence, we get the following properties.
\begin{proposition}\label{prop:nonexp}
	The operator $T$ is non-expansive, i.e. for any two functions $\mathbf{Q}_i:\overline{\cD}_N \to \R$, $i=1,2$, it holds
	\begin{equation*}
		\Norm{T(\mathbf{Q}_1)-T(\mathbf{Q}_2)}{\infty}\le \Norm{\mathbf{Q}_1-\mathbf{Q}_2}{\infty}.
	\end{equation*}
	Moreover, the operator $h$ is Lipschitz continuous with Lipschitz constant $2$.
\end{proposition}
\begin{proof}
	Being $h(\mathbf{Q})=T(\mathbf{Q})-\mathbf{Q}$, if we prove the first part of the statement, the second follows directly. To prove the first part, just observe that for any $(n,x,t,a) \in \mathbf{D}_N$ it holds
	\begin{align*}
		|T(\mathbf{Q}_1)(n,x,t,a)&-T(\mathbf{Q}_2)(n,x,t,a)|\\&\le \E_{\pi^{(a,n)},n,(x,t)}\left[\left|\max_{b \in D(X^{m}_{n+1},\gamma^{m}_{n+1})}\mathbf{Q}_1(n+1,X_{n+1}^{m},\gamma_{n+1}^{m},b)\right.\right.\\&\qquad \left.\left.-\max_{b \in D(X^{m}_{n+1},\gamma^{m}_{n+1})}\mathbf{Q}_2(n+1,X_{n+1}^{m},\gamma_{n+1}^{m},b)\right|\right]\\
		&\le \E_{\pi^{(a,n)},n,(x,t)}\left[\max_{b \in D(X^{m}_{n+1},\gamma^{m}_{n+1})}\left|\vphantom{\max_{b \in D(X^{m}_{n+1},\gamma^{m}_{n+1})}}\mathbf{Q}_1(n+1,X_{n+1}^{m},\gamma_{n+1}^{m},b)\right.\right.\\&\qquad \left.\left.\vphantom{\max_{b \in D(X^{m}_{n+1},\gamma^{m}_{n+1})}}-\mathbf{Q}_2(n+1,X_{n+1}^{m},\gamma_{n+1}^{m},b)\right|\right]\le \Norm{\mathbf{Q}_1-\mathbf{Q}_2}{\infty}, 
	\end{align*}
	where in the second inequality we used Lemma \ref{lemma:max}. \qed 
\end{proof}
Now let us define, for any $\beta>0$, the functions $h_\beta(\mathbf{Q})=\frac{h(\beta\mathbf{Q})}{\beta}$. For such a family of functions we have the following result.
\begin{proposition}
	Consider the operator $h_\infty:\overline{F}_N \to \overline{F}_N$ acting on any function $\mathbf{Q}:\overline{\cD}_N \to \R$ as
	\begin{equation*}
		h_\infty(\mathbf{Q})=h(\mathbf{Q})-\mathbf{r}.
	\end{equation*}
	Then $\lim_{\beta \to +\infty}h_\beta(\mathbf{Q})=h_\infty(\mathbf{Q})$ and the limit is uniform in compact sets of $\R^{|\overline{\cD}_N|}$.
\end{proposition}
\begin{proof}
	First of all, let us observe that $h_\infty(\mathbf{Q})$ is positively homogeneous by definition, i.e. $h_\infty(\beta \mathbf{Q})=\beta h_\infty(\mathbf{Q})$ for any $\beta>0$. Moreover, it holds $h(\mathbf{Q})=h_\infty(\mathbf{Q})+\mathbf{r}$, thus
	\begin{equation*}
		h_\beta(\mathbf{Q})=h_\infty(\mathbf{Q})+\frac{\mathbf{r}}{\beta}.
	\end{equation*}
	Hence it is clear that $\lim_{\beta \to +\infty}h_\beta(\mathbf{Q})=h_\infty(\mathbf{Q})$. To prove that the convergence is uniform in compact sets of $\R^{|\overline{\cD}_N|}$, let us first observe that since $h$ is Lipscthiz continuous with constant $2$, so it is $h_\infty$ and then also $h_\beta$. Moreover, if we assume $\beta>1$ it holds
	\begin{equation*}
		\Norm{h_\beta(\mathbf{Q})}{\infty}\le \Norm{h_\infty(\mathbf{Q})}{\infty}+\Norm{\mathbf{r}}{\infty}\le 2\Norm{\mathbf{Q}}{\infty}+\Norm{\mathbf{r}}{\infty},
	\end{equation*}
	where we used the fact that $h_\infty$ is $2$-Lipschitz and $h_\infty(0)=0$. Without loss of generality, we can consider just the family of compact sets given by $B_R=\{\mathbf{Q} \in \R^{|\overline{\cD}_N|}: \ \Norm{\mathbf{Q}}{\infty}\le R\}$. On $B_R$ we have that $h_\beta$ are equi-Lipschitz and equi-bounded, thus, in this case, pointwise convergence implies uniform convergence, concluding the proof. \qed 
\end{proof}
Let us observe that $h_\infty$ coincides with $h$ if we assume that the rewards are identically $0$. Now let us consider the following ODEs on $\R^{|\overline{\cD}_N|}$:
\begin{equation}\label{eq:ODE1}
	\der{\mathbf{Q}(t)}{t}=h(\mathbf{Q}(t)), \ t > 0
\end{equation}
and
\begin{equation}\label{eq:ODE2}
	\der{\mathbf{Q}(t)}{t}=h_\infty(\mathbf{Q}(t)), \ t > 0.
\end{equation}
Being $h$ and $h_\infty$ globally Lipschitz, these equations always admit a unique global solution for any initial data. With a similar strategy to the one adopted in \cite{abounadi2001learning}, we can prove the following stability result.
\begin{proposition}\label{prop:ODE}
	Equation \eqref{eq:ODE1} admits $\mathbf{Q}^*$ as its unique globally asymptotically stable equilibrium. Equation \eqref{eq:ODE2} admits $0$ as its unique globally asymptotically stable equilibrium.
\end{proposition}
\begin{proof}
	Clearly, Equation \eqref{eq:ODE2} is a particular case of \eqref{eq:ODE1}, with null rewards and $\mathbf{Q}^* =0$. Thus, if we prove the first statement, we obtain the second one as a direct consequence. Let us stress that the equation $h(\mathbf{Q})=0$ is actually a recasting of Bellman's Equation \eqref{eq:BEqn}, thus it is clear that $\mathbf{Q}^*$ is the unique equilibrium of the system. To prove that $\mathbf{Q}^*$ is globally asymptotically stable, observe that we can rewrite equation \eqref{eq:ODE1} as
	\begin{equation*}
		\der{\mathbf{Q}(t)}{t}=T(\mathbf{Q}(t))-\mathbf{Q}(t), \ t > 0,
	\end{equation*}
	where we proved in Proposition \ref{prop:nonexp} that $T$ is non-expansive. Thus, by \cite[Theorem $3.1$]{borkar1997analog}, we know that any solution $\mathbf{Q}(t)$ converges towards an equilibrium point of the ODE (possibly depending on the initial data). However, the ODE admits a unique equilibrium point $\mathbf{Q}^*$, hence we get, for any initial data, $\lim_{t \to \infty}\Norm{\mathbf{Q}(t)-\mathbf{Q}^*}{\infty}=0$.\qed 
\end{proof}
Now we move to the study of the properties of $\mathbf{M}^m$. Let us first define the filtration $\overline{\cF}_m=\sigma(\mathbf{Q}^i,\mathbf{M}^i, \ i \le m)$. We want to prove, as done in \cite[Proposition $3$, item $1$]{bhatnagar2021finite}, the following property.
\begin{proposition}\label{martingale}
	Suppose $\mathbf{Q}_0$ is square-integrable. Then the process $(\mathbf{M}^m)_{m \in \N}$ is a square-integrable martingale difference sequence with respect to the filtration $\overline{\cF}_m$. Moreover, for any initial condition $\mathbf{Q}_0$, it holds
	\begin{equation}\label{eq:mart}
		\E[\Norm{\mathbf{M}^{m+1}}{\infty}^2|\overline{\cF}_m]\le 4\Norm{\mathbf{Q}^m}{\infty}^2 \quad \forall m \in \N_0
	\end{equation}
	almost surely.
\end{proposition}
\begin{proof}
	It is clear by definition that $\mathbf{M}^m$ is $\overline{\cF}_m$-measurable for any $m \in \N$.  Let us prove that, for any $m \in \N$, $\E[\mathbf{M}^{m+1}|\overline{\cF}_m]=0$. Indeed since $\mathbf{Q}^m$ is $\overline{\cF}_m$-measurable while $(X^m,\gamma^m)$ is independent of $\overline{\cF}_m$ we get, for any $(n,x,t,a) \in \overline{\cD}_N$ with $n \le N-1$,
	\begin{align*}
		&\E[\mathbf{M}^{m+1}(n,x,t,a)|\overline{\cF}_m]\\&=\E\left[\left.\max_{b \in D(X^{m}_{n+1},\gamma^{m}_{n+1})}\mathbf{Q}^m(n+1,X_{n+1}^{m},\gamma_{n+1}^{m},b)\right|\overline{\cF}_m\right]\\
		&-\E\left[\left.\E_{\pi^{(a,n)},n,(x,t)}\left[\left.\max_{b \in D(X^{m}_{n+1},\gamma^{m}_{n+1})}\mathbf{Q}^m(n+1,X_{n+1}^{m},\gamma_{n+1}^{m},b)\right| \mathbf{Q}^m\right]\right| \overline{\cF}_m\right]=0,
	\end{align*}
	while the equality is trivial if $n=N$.\\
	Now, we need to prove that $\mathbf{M}^{m}$ is square-integrable. To do this, observe that, for any $(n,x,t,a) \in \overline{\cD}_N$ with $n \le N-1$ and any $m \in \N_0$, we have, according to \eqref{Qlearnfunctional},
	\begin{align*}
		|\mathbf{Q}^{m+1}(n,x,t,a)|&\le (1-\alpha(m))\Norm{\mathbf{Q}^m}{\infty}+\alpha(m)\Norm{\mathbf{Q}^m}{\infty}+\alpha(m)\Norm{\mathbf{r}}{\infty}\\& \le \Norm{\mathbf{Q}^m}{\infty}+\Norm{\mathbf{r}}{\infty},
	\end{align*}
	and the same inequality clearly holds even if $n=N$. Thus we get
	\begin{equation*}
		\Norm{\mathbf{Q}^{m+1}}{\infty}\le \Norm{\mathbf{Q}^m}{\infty}+\Norm{\mathbf{r}}{\infty}.
	\end{equation*}
	Iterating the previous relation $m$ times and taking the square we get
	\begin{equation*}
		\Norm{\mathbf{Q}^{m}}{\infty}^2\le 2(\Norm{\mathbf{Q}^0}{\infty}^2+m^2\Norm{\mathbf{r}}{\infty}^2).
	\end{equation*}
	Hence, by definition of $\mathbf{M}^{m+1}$ as in \eqref{eq:M} we have, for any $(n,x,t,a) \in \overline{\cD}_N$,
	\begin{equation}\label{premart}
		|\mathbf{M}^{m+1}(n,x,t,a)|^2\le 4\Norm{\mathbf{Q}^m}{\infty}^2\le 8(\Norm{\mathbf{Q}^0}{\infty}^2+m^2\Norm{\mathbf{r}}{\infty}^2).
	\end{equation}
	Taking the maximum over $(n,x,t,a) \in \overline{\cD}_N$ and then the expectation we conclude that
	\begin{equation*}
		\E[\Norm{\mathbf{M}^{m+1}}{\infty}^2]\le 8(\E[\Norm{\mathbf{Q}^0}{\infty}^2]+m^2\Norm{\mathbf{r}}{\infty}^2)<\infty.
	\end{equation*}
	Moreover, if we take the maximum over $(n,x,t,a) \in \overline{\cD}_N$ and then the conditional expectation with respect to $\overline{\cF}_m$ in inequality \eqref{premart}, we achieve inequality \eqref{eq:mart}. \qed 
\end{proof}
Now we are ready to prove the convergence of the $Q$-learning algorithm under the \textit{tapering stepsize} assumption, as done in \cite{bhatnagar2021finite}.
\begin{proposition}\label{prop:conv}
	Let $\alpha:\N_0 \to \R$ be such that $\sum_{m=0}^{\infty}\alpha(m)=\infty$ and $\sum_{m=0}^{\infty}\alpha^2(m)<\infty$. Then $\mathbf{Q}^{m} \to \mathbf{Q}^*$ almost surely, i.e.
	\begin{equation*}
		\bP(\lim_{m \to \infty}\mathbf{e}(m)=0)=1.
	\end{equation*} 
\end{proposition}
\begin{proof}
	By Proposition \ref{prop:nonexp} we know that $h$ is Lipschitz continuous, while Proposition \ref{prop:ODE} tells us that the ODE \eqref{eq:ODE2} admits the origin as its globally asymptotically stable equilibrium. Moreover, Proposition \ref{martingale} ensures that $(\mathbf{M}^m)$ is a square-integrable martingale difference sequence satisfying inequality \eqref{eq:mart}. Thus, the claim follows by \cite[Theorem $2.2$]{borkar2000ode}. \qed
\end{proof}
On the other hand, we could also ask what happens in the \textit{constant stepsize} assumption. This is stated in the following proposition, whose proof is analogous to the previous one, while resorting to \cite[Theorem $2.3$]{borkar2000ode}.
\begin{proposition}\label{prop:conv1}
	Let $\alpha:\N_0 \to \R$ be constant. Then for any $\varepsilon>0$ there exists a constant $b_1(\varepsilon)$ independent of $\alpha$ such that
	\begin{equation*}
		\limsup_{n \to \infty}\bP(\mathbf{e}(m)\ge \varepsilon)\le b_1(\varepsilon)\alpha.
	\end{equation*}
	If, furthermore, $\mathbf{Q}^*$ is globally exponentially asymptotically stable, then there exists a constant $b_2$ independent of $\alpha$ such that
	\begin{equation*}
		\limsup_{n \to \infty}\E[\mathbf{e}^2(m)]\le b_2\alpha.
	\end{equation*}
\end{proposition}
Clearly, we cannot guarantee convergence if $\alpha$ is constant, but we can \textit{modulate} in some sense the absolute error in such a way that the result is a \textit{good enough} approximation. Different choices of $\alpha(m)$ could also affect the \textit{speed of convergence}, as we will see in the following Section, thus a first study on a constant learning rate could be useful for tuning the right form of a tapering one. The balance between absolute error and speed of convergence is, for instance, expressed in \cite[Theorem 2.4]{borkar2000ode}, but relies on further properties of the sequence $(\mathbf{Q}^m)_{m \in \N}$, that will be investigated in future works. In the next section we will provide two toy examples to our theory and we will give some simulation results, showing the effects of different choices of the learning rate on the average total reward after a certain number of iterations.
\section{Numerical tests}\label{sec5}
\subsection{Switching uneven coins}\label{sec51}
\subsubsection{Theoretical setting}\label{sec511}
Imagine we want to play a game of heads or tails against an artificial intelligence (AI), in which we can only bet on heads while the AI can only bet on tails. From the point of view of the AI, it gains a point each time it is tails, otherwise it loses a point. However, it can use a certain number of different uneven coins and can \textit{switch} between them at any turn of the game. Moreover, the opponent is human and its conception of randomness is not perfect, but subject to \textit{gambler's fallacy} (see, for instance, \cite{warren2018re} and references therein). Thus, if it is tails for \textit{too many} subsequent turns, he will start thinking that the AI is cheating. In a certain sense, the AI is strongly penalized if too many subsequent tails come out. On the other hand, this is not true for the heads, as the human opponent does not expect the AI to \textit{cheat against itself}.\\
From a mathematical point of view, if we consider $X_n$ to be the result of the $n$-th coin toss and $\gamma_n$ to be the sojourn time process of $X_n$, the AI receives a reward (or penalization) depending on both $X_n$ and $\gamma_n$ via the function
\begin{equation*}
	r(x,t)=\begin{cases} 1 & x=0, \ t\le t_{\rm cheat}\\
		-1 & x=1 \\
		r_{\rm cheat} & x=0, \ t>t_{\rm cheat}
	\end{cases}
\end{equation*}
where $t_{\rm cheat}$ is the maximum number of subsequent tails before the human opponent thinks the AI is cheating and $r_{\rm cheat}<0$ is the penalization that the AI receives if the human suspects it of cheating. Suppose we have a number $M$ of coins identified with the actions $A=\{a_1,\dots,a_{M}\}$. For each coin $a_i$, let us denote by $p_i$ the probability of success (i.e. the probability that, flipping the coin $a_i$, it is heads). Let us also identify tails with $0$ and heads with $1$. Then we define the transition kernel $P$ as
\begin{align*}
	P(1,0|0,t,a_i)=p_i && P(0,t+1|0,t,a_i)=1-p_i \\
	P(1,t+1|1,t,a_i)=p_i && P(0,0|1,t,a_i)=1-p_i
\end{align*}
for any $i=1,\dots,M$. Let us also fix the horizon (i.e. the number of coin tossing for each play) $N>0$.\\
Let us stress out that since $r$ depends on the sojourn time, we should expect, in general, that the optimal policy could also depend on it. This means that we cannot presume, in general, $X=\{X_n, \ n \in \N\}$ to be a Markov process, as shown in the following result.
\begin{proposition}
	For the \emph{switching coins problem} with $N>1$, the following statements are equivalent:
	\begin{enumerate}
		\item $X$ is $\bP_{\pi,(x_0,0)}$-Markov for any policy $\pi \in F^N$ and any initial state $x_0=0,1$.
		\item It holds $p_1=p_2=\dots=p_M$.
	\end{enumerate}
\end{proposition}
\begin{proof}
	The implication $(2) \Rightarrow (1)$ is clear as the action has no real effect on the model, that becomes a classical Bernoulli scheme.\\
	Let us show that $(1) \Rightarrow (2)$. To do this, fix any $i,j \le M$ with $i \not = j$. Let us consider any policy $\pi$ such that $\pi_0(0,0)=a_i$ and $\pi_1(0,1)=a_j$. Let $T$ be the first inter-jump time of $X$. Being $X$ a Markov process w.r.t $\bP_{\pi,(0,0)}$, $T$ must be a geometric random variable. In general, it must hold
	\begin{equation*}
		\bP_{\pi,(0,0)}(T=1)=\bP_{\pi,(0,0)}(X_1=1,\gamma_1=0|X_0=0,\gamma_0=0)=p_i,
	\end{equation*}
	and
	\begin{equation*}
		\bP_{\pi,(0,0)}(T=2)=\bP_{\pi,(0,0)}(X_2=1,\gamma_2=0|X_1=0,\gamma_1=1)\bP(T>1)=p_j(1-p_i).
	\end{equation*}
	However, being $T$ a geometric random variable, there exists $\overline{p} \in [0,1]$ such that
	\begin{equation*}
		\begin{cases}
			\overline{p}=p_i\\
			(1-\overline{p})\overline{p}=(1-p_i)p_j,
		\end{cases}
	\end{equation*}
	that implies either $p_i=p_j$ or $p_i=1$. Exchanging the roles of $i$ and $j$ we conclude that, in any case, $p_i=p_j$. Finally, being $i,j$ arbitrary, we conclude the proof. \qed 
\end{proof}
\subsubsection{Numerical results}\label{sec512}
To exhibit some features of the method that we have presented here, we consider the specific case of $M=2$ coins, with parameters $p_1=1/5$ and $p_2=4/5$, $t_{\rm cheat}=3$ and $r_{\rm cheat}=-10$. Moreover, we set the horizon as $N=200$. To obtain the optimal policy, we adopt the Q-learning algorithm based on the recursion formula \eqref{Qlearn}. At the $(m+1)$-th iteration, for the $n$-th turn, the action $a$ on which the evaluation is taken in $Q_n^{m+1}(x,t,a)$ is selected in the set $\arg\max Q_{n+1}^{m}(X_{n+1}^m,\gamma_{n+1}^m,\cdot)$, while $(x,t)=(X_{n}^m,\gamma_n^m)$, where both $(X_{n+1}^m,\gamma_{n+1}^m)$ and $(X_{n+1}^m,\gamma_{n+1}^m)$ are obtained via stochastic simulation. All the other values are left unchanged. This procedure is repeated for $m=0,\dots,N_{\rm eps}-1$ where $N_{\rm eps}=2000$ is the number of episodes. At $m=0$, the function $Q^0$ is initialized randomly. Moreover, in the following figures the learning rate $\alpha$ in Equation \eqref{Qlearn} will be denoted as ${\rm LR}$. \\
The total reward is evaluated on each simulated trajectory $(X_n^m,\gamma_n^m,a_n^m)$ (where $a_n^m$ is the action taken at step $n$ in episode $m$) as follows:
\begin{equation*}
	R_m=\sum_{n=0}^{N}r(X_n^m,\gamma_n^m)
\end{equation*}
so that the expected reward can be approximated via Monte Carlo methods. Here we consider $100$ different evaluations and then we take the average. To smooth further the data, we do not consider directly $\E[R_m]$ but the expected value of the average on batches of $50$ episodes, i.e. the \textit{average expected reward} defined as 
\begin{equation*}
	R_{k,{\rm avg}}=\E\left[\frac{1}{50}\sum_{m=50k}^{50(k+1)-1} R_m\right].
\end{equation*}
The results for different values of the learning rate can be seen in Figure \ref{fig:figure1}, while in Figure \ref{fig:figure2} we plot only the trajectory of $R_{k,{\rm avg}}$ for ${\rm LR}=0.3$. As we can see on both figures (and as expected), the average expected reward is increasing with $m$. Let us recall that, by Theorem \ref{thm:conv}, we do not presume that the method converges, but it still provides a quite good approximation of the solution. As one can see from the plots, we obtain a lower average expected reward for ${\rm LR}=0.1,0.2$ while the values for ${\rm LR}=0.3,0.4$ and $0.5$ are higher and quite close each other. This seems to contradict Theorem \ref{thm:conv}: actually, this phenomenon is due to trade-off between accuracy and speed of convergence (that could be explained, for instance, by \cite[Theorem 2.4]{borkar2000ode}). Thus, even if for ${\rm LR}=0.1$ the method should be more accurate, it takes also much more episodes to reach such an accuracy. For this reason, we use ${\rm LR}=0.3$ as reference case.
\begin{center}
	\begin{minipage}{0.49\linewidth}
		\includegraphics[width=\linewidth]{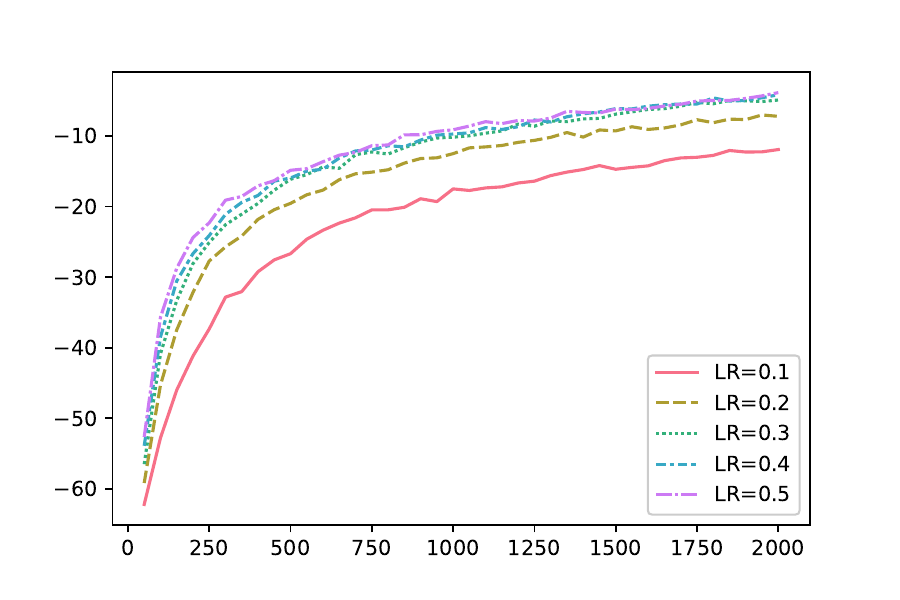}
		\captionof{figure}{Plot of the average expected rewards for different values of the learning rate in the \textit{switching coins} example.}
		\label{fig:figure1}
	\end{minipage}
	\begin{minipage}{0.49\linewidth}
		\includegraphics[width=\linewidth]{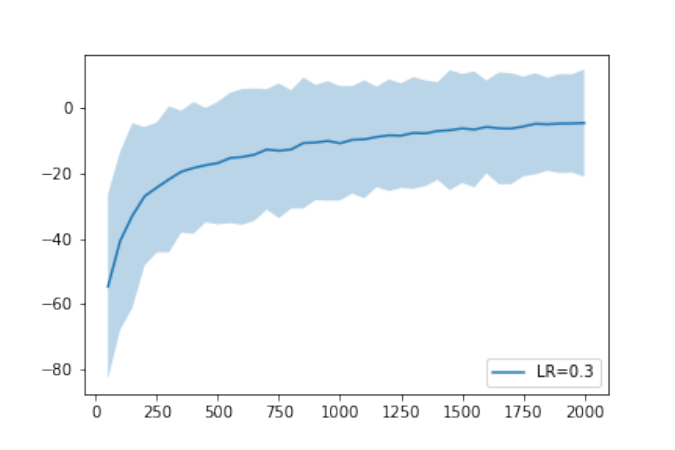}
		\captionof{figure}{Plot of the average expected rewards for ${\rm LR}=0.3$ with a $95\%$ confidence interval in the \textit{switching coins} example.}
		\label{fig:figure2}
	\end{minipage}
\end{center}
Together with the average, we considered also the \textit{minimal} and the \textit{maximal expected rewards} defined respectively by
\begin{equation*}
	R_{k,{\rm min}}=\E[\min_{50k\le m < 50(k+1)}R_m] \qquad \mbox{ and }\qquad  R_{k,{\rm max}}=\E[\max_{50k\le m < 50(k+1)}R_m].
\end{equation*}
Again, both of them are obtained via Monte Carlo methods on $100$ evaluations. The plots of the minimal expected rewards are presented in Figure \ref{fig:figure3} and the relative one for ${\rm LR}=0.3$ is given in Figure \ref{fig:figure4}. It is interesting to observe the higher variability of the trajectory for ${\rm LR}=0.3$ with respect to the one for the average expected reward. As expected, all the trajectories are increasing and remain negative (as the average ones are also negative). Moreover, even in this case we observe that the values for ${\rm LR}=0.1$ and $0.2$ are lower than the ones for ${\rm LR}=0.3,0.4$ and $0.5$, the latter being still quite close each other. This justifies the choice of ${\rm LR}=0.3$ as a reference value even if the \textit{smoothing procedure} is made by \textit{taking the minimum} instead of \textit{averaging} on batches of subsequent episodes.
\begin{center}
	\begin{minipage}{0.49\linewidth}
		\includegraphics[width=\linewidth]{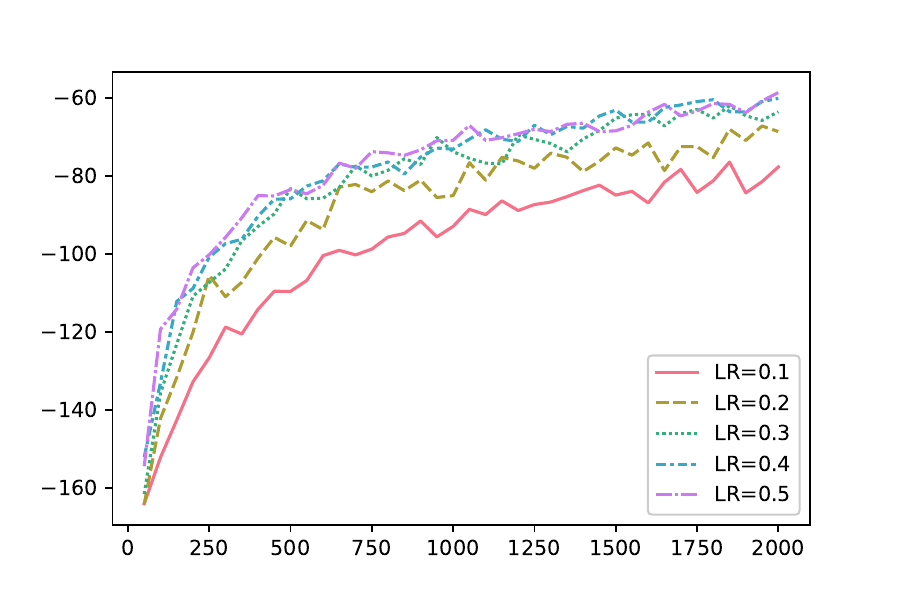}
		\captionof{figure}{Plot of the minimal expected rewards for different values of the learning rate in the \textit{switching coins} example.}
		\label{fig:figure3}
	\end{minipage}
	\begin{minipage}{0.49\linewidth}
		\includegraphics[width=\linewidth]{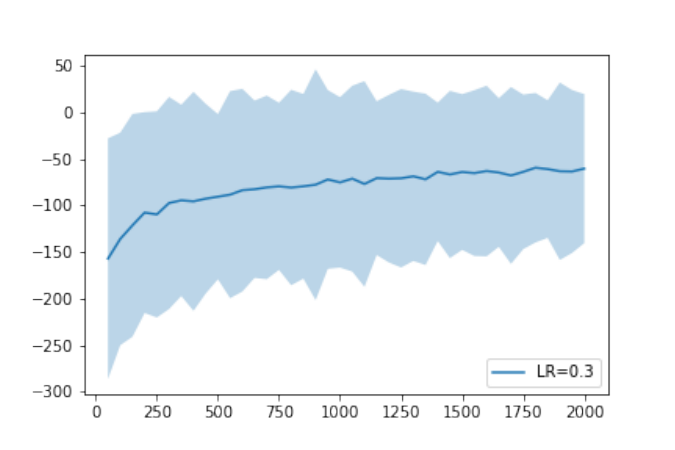}
		\captionof{figure}{Plot of the minimal expected rewards for ${\rm LR}=0.3$ with a $95\%$ confidence interval in the \textit{switching coins} example.}
		\label{fig:figure4}
	\end{minipage}
\end{center}
Similar observations can be carried on for the maximal expected rewards, depicted in Figures \ref{fig:figure5} and \ref{fig:figure6}. This time, the only evidently lower value is given by ${\rm LR}=0.1$. Again, we have higher variability with respect to the average one. Moreover, we deduce from the graphs that, in this specific case, the maximal expected rewards are non-negative already for $m>250$, which is an information we could not deduce from the average ones. 
\begin{center}
	\begin{minipage}{0.49\linewidth}
		\includegraphics[width=\linewidth]{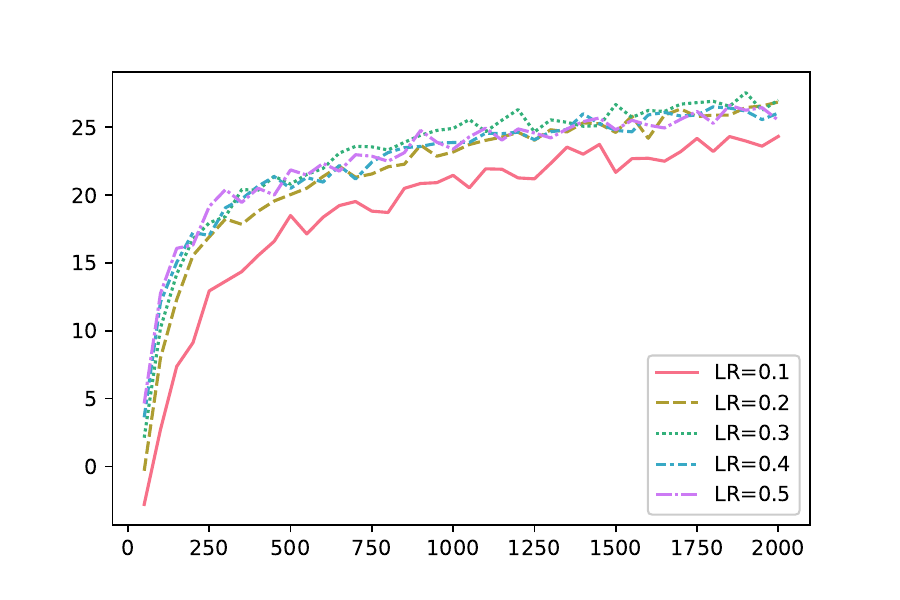}
		\captionof{figure}{Plot of the maximal expected rewards for different values of the learning rate in the \textit{switching coins} example.}
		\label{fig:figure5}
	\end{minipage}
	\begin{minipage}{0.49\linewidth}
		\includegraphics[width=\linewidth]{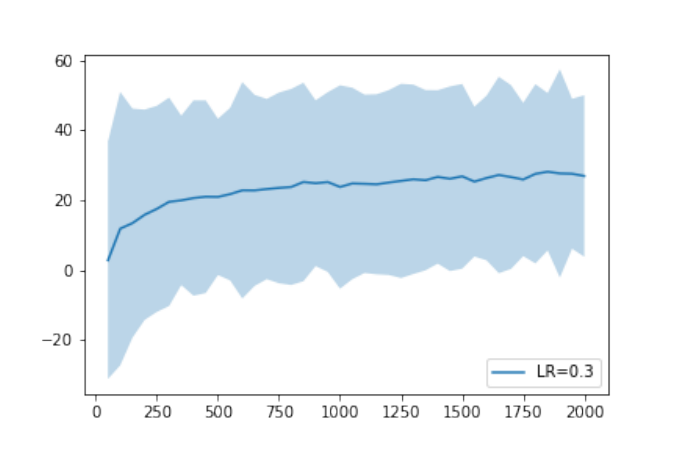}
		\captionof{figure}{Plot of the maximal expected rewards for ${\rm LR}=0.3$ with a $95\%$ confidence interval in the \textit{switching coins} example.}
		\label{fig:figure6}
	\end{minipage}
\end{center}
To ensure convergence, one could consider a variable learning rate. It is clear that, to guarantee a certain \textit{speed of convergence}, one can construct an \textit{ad-hoc} variable learning rate according to the behaviour of the algorithm in the constant learning rate case. Precisely, let us set
\begin{equation*}
	{\rm LR}(m)=\frac{1}{\lceil \frac{m+1}{100}\rceil +1}
\end{equation*}
so that we can use Proposition \ref{prop:conv} to guarantee convergence. Again, we compare this case with our reference case ${\rm LR}=0.3$. In particular one has that ${\rm LR}(m)<0.3$ as soon as $\lceil \frac{m+1}{100}\rceil \ge 3$, that is to say $m \ge 200$. This can be seen explicitly in Figures \ref{fig:figure7}, \ref{fig:figure8} and \ref{fig:figure9} respectively for the average, minimal and maximal expected rewards. Indeed, in all the cases, despite the expected reward is increasing as we desire, the variable learning rate case is faster than the fixed one at the start, while it becomes slower around $m=200$. However, the plots suggest that different choices for ${\rm LR}(m)$ could lead to a fast and, at the same time, accurate solution.\\
\begin{center}
	\begin{minipage}{0.49\linewidth}
		\includegraphics[width=\linewidth]{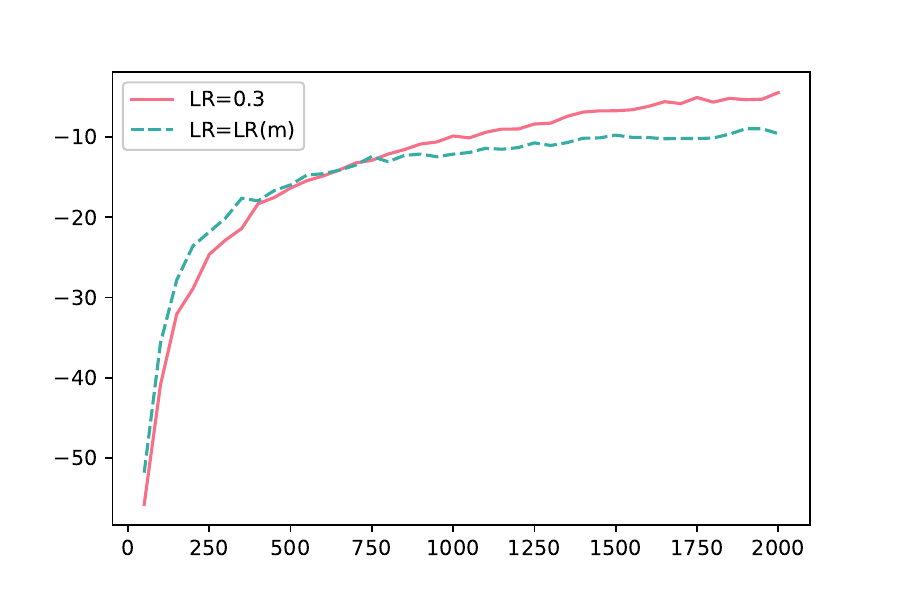}
		\captionof{figure}{Plot of the average expected rewards for fixed and variable learning rates in the \textit{switching coins} example.}
		\label{fig:figure7}
	\end{minipage}
	\begin{minipage}{0.49\linewidth}
		\includegraphics[width=\linewidth]{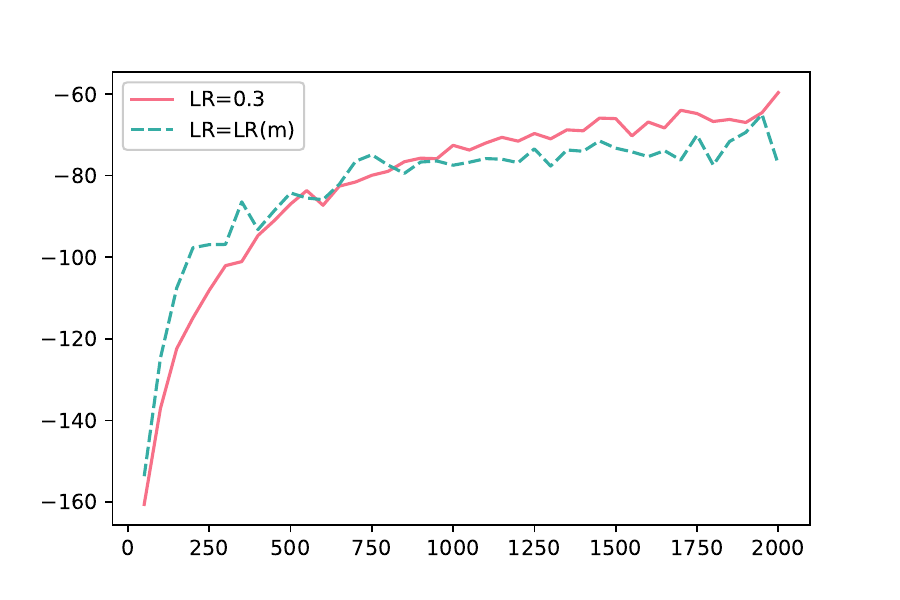}
		\captionof{figure}{Plot of the minimal expected rewards for fixed and variable learning rates in the \textit{switching coins} example.}
		\label{fig:figure8}
	\end{minipage}
\end{center}
\begin{center}
	\begin{minipage}{0.49\linewidth}
		\includegraphics[width=\linewidth]{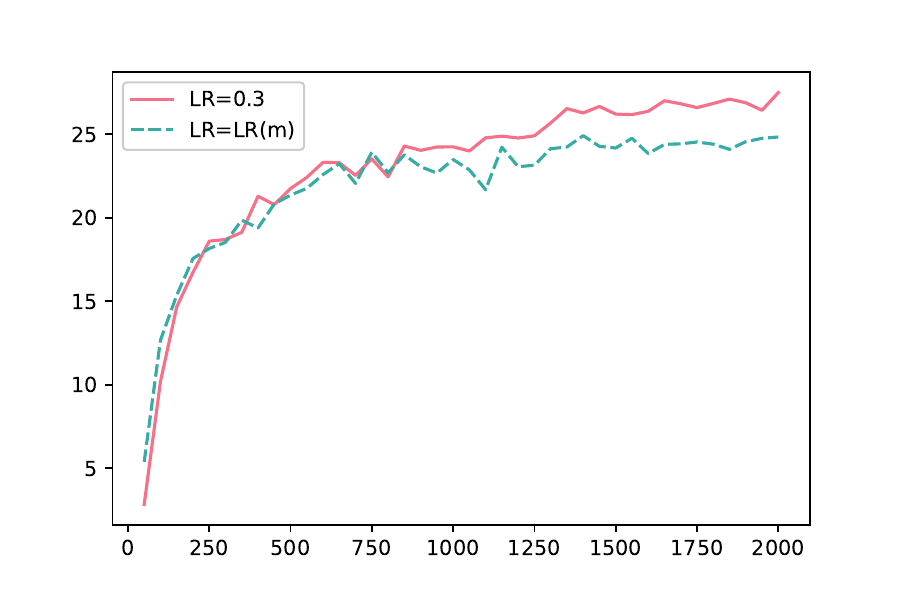}
		\captionof{figure}{Plot of the maximal expected rewards for fixed and variable learning rates in the \textit{switching coins} example.}
		\label{fig:figure9}
	\end{minipage}
	\begin{minipage}{0.49\linewidth}
		\includegraphics[width=\linewidth]{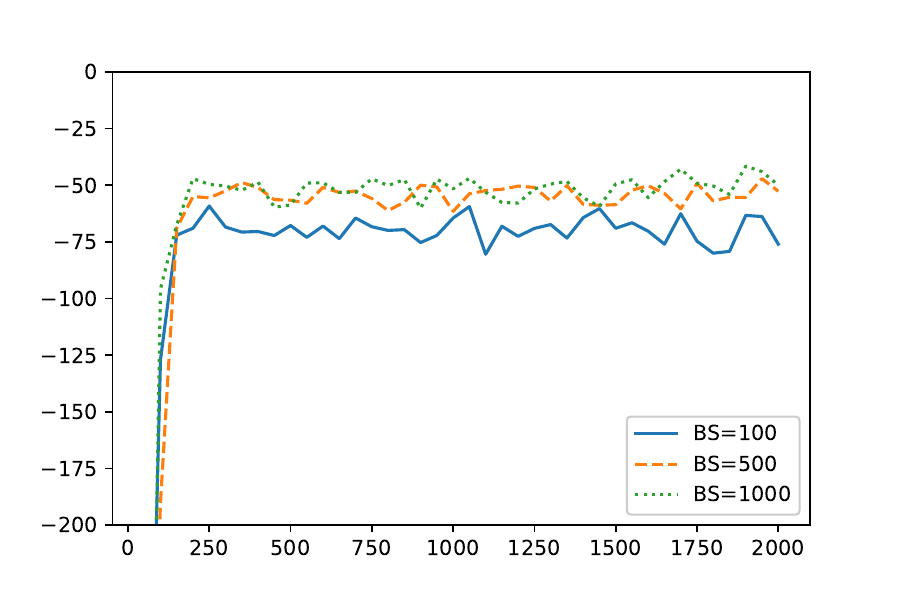}
		\captionof{figure}{Plot of the average rewards obtained via $Q$-network for ${\rm NN}=128$, ${\rm NL}=3$ and different values of ${\rm BS}$ in the \textit{switching coins} example.}
		\label{fig:figure10}
	\end{minipage}
\end{center}
Let us stress that even if we want to save the value of the function $Q$ only in the \textit{reachable states}, the dimension of the $Q$-table becomes quite difficult to work with. Indeed, denoting by $C(|E|,N,|A|)=|\overline{E}_N|$ the cardinality of the set of reachable states, it holds
\begin{equation*}
	C(|E|,N,|A|)=\frac{|E||A|N(N+1)}{2},
\end{equation*}
so, despite being linear in $|E|$ and $|A|$, it is quadratic in $N$. This increases drastically the space complexity of the method with respect to the horizon.\\
A possible solution for such a space complexity problem could be the implementation of a different function approximation method for $Q$. Precisely, let us first observe that we are actually trying to obtain an approximation of the function $\mathbf{Q}^*$ defined on $\overline{\cD}_N$. In place of using the $Q$-learning algorithm based on Equation \eqref{Qlearn}, we can implement a neural network to provide such a function approximation. In this way, the (almost) only things we have to store in memory are the weights of the network. The input layer of the network is constituted by three neurons (since the input for a fixed action $a$ is given by $(n,x,t)$), while the output layer is given by a number of neurons equal to the number of actions $|A|$. Indeed, the $i$-th output neuron will give a normalized approximation of $\mathbf{Q}(n,x,t,a_i)$, which will be used to determine the optimal action.\\
Here, we add a number ${\rm NHL}$ of Dense hidden layers, each one composed of ${\rm NN}$ neurons. The activation function of the hidden layers is a Rectified Linear Unit, while for the output layers we adopt a Softmax function. As loss function, we consider the Mean Squared Error one, while the optimization procedure is provided by the Adam procedure (see \cite{kingma2014adam}). As usual (see, for instance, \cite{mnih2015human}), to train the network we need to store some trajectories (with their rewards and policies) in a replay memory. We fix its dimension as $N_{\rm RM}=5000$, while we start the training session as soon as the replay memory reaches the length $N_{\rm MRM}=4000$. The training is carried on a batch of instances in the replay memory of size ${\rm BS}$ that are randomly selected from it. Moreover, let us stress that two actual networks are used: one that is fitted on each step (the actual \textit{model}), while a second one (called the \textit{target model}) is updated only each $100$ fits of the first one and it is used for prediction during the training sessions of the \textit{model}. This procedure is shown to reduce the instability of the non-linear function approximation (see \cite{mnih2015human}). Obviously, as a trade-off for the space complexity, the $Q$-network approximation takes much more time than the direct $Q$-learning algorithm, both in terms of convergence (the $Q$-network approximation exhibits a slower speed of convergence to the optimal value) than in terms of the execution time (at each step, the $Q$-network has to fit on a training set of ${\rm BS}$ states: this takes a longer time to compute even on the same number $N_{\rm eps}$ of episodes). Moreover, a careful study on the parameters of the network has to be carried on to ensure its correct functioning.\\
Let us first fix ${\rm NHL}=3$ and ${\rm NN}=128$ and let us consider different values of ${\rm BS}$. Since the horizon is fixed to $N=200$, we expect the algorithm to work better as ${\rm BS} \ge 200$. Indeed, the states that are stored in the replay memory correspond not only to different couples $(x,t)$, but also to different steps $n$ (as the function we want to approximate depends also on $n$). If we consider ${\rm BS}<200$, then, even if we do not carry a random selection in the replay memory, we cannot have an instance for each step of the process. In general, we expect the function approximation to work better as the batch size ${\rm BS}$ is taken to be bigger than the horizon $N$. In Figure \ref{fig:figure10} we plotted the average rewards $\widetilde{R}_{k,\rm avg}=\frac{1}{50}\sum_{m=50k}^{50(k+1)-1}R_m$ (we evaluated it on a single simulation in place of taking the expected value via Monte Carlo methods) obtained in this way for three different values of ${\rm BS}$. Here we see that, as we supposed, the case in which ${\rm BS}=100$ provides a worst reward than the other two. Moreover, the cases with ${\rm BS}=500$ and $1000$ provide a quite similar reward (with ${\rm BS}=1000$ being only slightly better). Thus, for the subsequent tests, we consider ${\rm BS}=500$. Moreover, the average reward seems to stabilize after the first episodes, thus we reduce, from now on, the number of episodes to $N_{\rm eps}=1000$.\\
Now, we fix ${\rm BS}=500$ and ${\rm NHL}=3$ and we consider different values for ${\rm NN}$. The plots of the average rewards are shown in Figure \ref{fig:figure11}. Here we can see that reducing the number of neurons per layer actually reduces the average reward in the first episodes, but it does not seem to have any big impact on the average reward after the first $400$ episodes. Since we stop training after $1000$ episodes, we presume that ${\rm NN}=8$ are enough neurons per layer.\\
Finally, we fix ${\rm BS}=500$ and ${\rm NN}=8$ and we consider different values for ${\rm NHL}$. The plots of the average rewards are given in Figure \ref{fig:figure12}. Here we see that, again, for low values of ${\rm NHL}$ (i.e. $2$ and $3$), the average reward is lower on the first episodes, but stabilizes on the same value of the case ${\rm NHL}=4$ after $400$ episodes circa. On the other hand, the \textit{non-monotonic} initial behaviour is avoided for ${\rm NHL}=4$ and the stabilization is achieved already around $300$ episodes. This suggests that, on one hand, it could be useful to set a lower number of hidden layers to obtain a certain \textit{gain on time}, while still considering an high number of episodes; on the other hand, one could also consider a higher number of hidden layers and then reduce the number of episodes on which the network has to be trained.
\begin{center}
	\begin{minipage}{0.49\linewidth}
		\includegraphics[width=\linewidth]{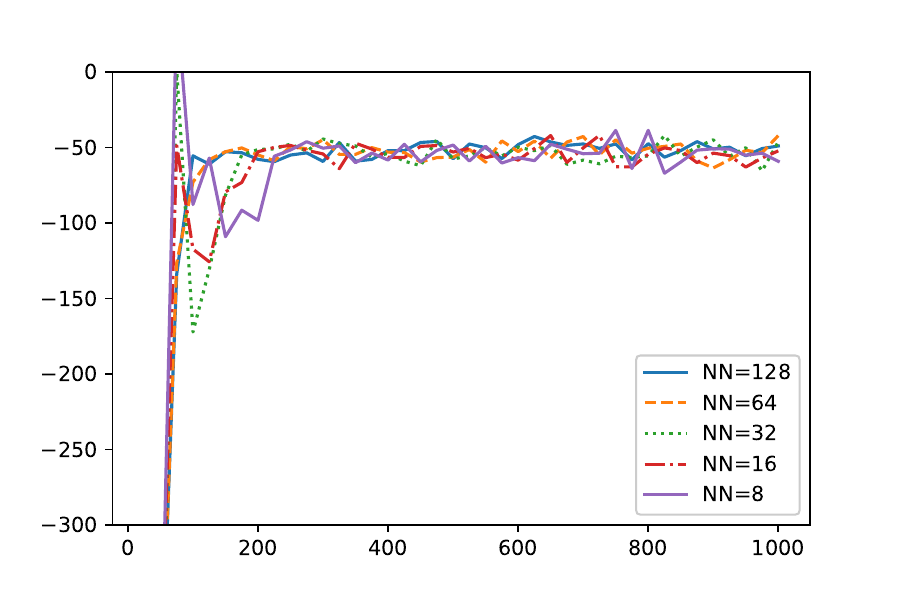}
		\captionof{figure}{Plot of the average rewards obtained via $Q$-network for ${\rm BS}=500$, ${\rm NHL}=3$ and different values of ${\rm NN}$ in the \textit{switching coins} example.}
		\label{fig:figure11}
	\end{minipage}
	\begin{minipage}{0.49\linewidth}
		\includegraphics[width=\linewidth]{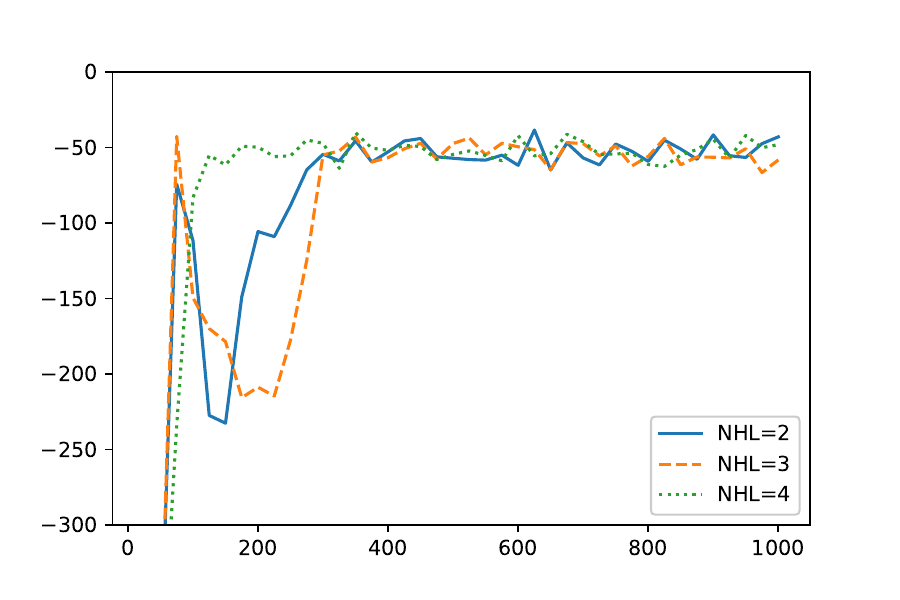}
		\captionof{figure}{Plot of the average rewards obtained via $Q$-network for ${\rm NN}=8$, ${\rm BS}=500$ and different values of ${\rm NHL}$ in the \textit{switching coins} example.}
		\label{fig:figure12}
	\end{minipage}
\end{center}
Clearly, the $Q$-network approach needs a deeper study. Up to now, it provides a lower expected reward with respect to the direct $Q$-learning method, despite this first light preliminary study on the hyper-parameters. We plan to improve the $Q$-network approach in future works.
\subsection{Preventive Maintenance with partially observable state}\label{sec52}
\subsubsection{Theoretical setting}\label{sec521}
Let us assume that we have a machine that can produce a good according to two different regimes $a_0$ and $a_1$. Moreover, suppose that we can \textit{observe} such a system at some (discrete) time instants $t_n=n$ and that it has to be shut down once we reach $t_N=N$. Let $Y_n$ be a process describing a certain \textit{state} of the system (for instance, the temperature of the machine). To give an idea of its role, we refer here to $Y_n$ as the \textit{stress level} of the system. We ask that $Y_0$ is fixed to a value that we call \textit{default stress level} and $Y_n \ge Y_0$ almost surely. We expect a more likely breakdown of the system if the \textit{stress level is too high}, i.e. when $Y_n$ exceeds a certain critical threshold. Thus, let us fix it as $\overline{Y}>Y_0$. Once $Y_n$ surpasses $\overline{Y}$, there is a certain probability that the system will break down, which increases with the time $Y_n$ spends above the default stress level. Precisely, let us suppose that as soon as $Y_n$ crosses $\overline{Y}$, if we denote by $\gamma_n$ the time spent above $Y_0$ right after this event up to time $n$, the machine breaks down as soon as $\gamma_n \ge \tau_1$, where $\tau_1$ is a non-negative integer-valued random variable with known distribution independent of $Y_n$, which we call \textit{stress threshold}. Once the system is broken, one needs to repair it and to reduce $Y_n$ to the default stress level $Y_0$. If we denote again by $\gamma_n$ the time that has passed since the last system broke down up to time $n$, we assume that the system is repaired as soon as $Y_n=Y_0$ and $\gamma_n \ge \tau_2$, where $\tau_2$ is a non-negative integer-valued random variable with known distribution independent of $Y_n$, assumed to be the \textit{time to repair}. Up to now, we did not link the behaviour of $Y_n$ with the two regimes $a_1$ and $a_2$. Let us assume that $Y_n$ is a Markov chain whose transition matrix depends on the choice of the regime $a_i$. The two regimes $a_1$ and $a_2$ can be considered as a \textit{slow} and \textit{fast} regimes. Under $a_1$ we have a certain profit-per-time-unit $c_1 \in \R$ and $Y_n$ is more likely to decrease or stay still. Under $a_2$, we have a different profit-per-time-unit $c_2>\max\{0,c_1\}$, however $Y_n$ is more likely to increase. If the system is broken, we have a cost that is proportional to the inactivity time of the system, represented via a negative profit-per-time-unit $c_3<\min\{c_1,-c_2,0\}$. This is a \textit{high-risk high-reward} system: using the \textit{fast} regime leads to a bigger profit than the \textit{slow} regime, while at the same time we are more exposed to possible breakdowns.\\
However, we further suppose that we cannot observe directly $Y_n$, but only a functional $(X_n)_{n \in \N}$ of it, taking values in $\{0,1,2\}$. To define $X_n$ we proceed as follows. First we fix $X_0=0$ and we denote by $\gamma_n$ its sojourn time. Suppose we already know the value of $(X_n,\gamma_n)$ and let us show how we determine $X_{n+1}$. If $X_n=0$, then we set
\begin{equation*}
	X_{n+1}=\begin{cases} 0 & Y_{n+1}\le \overline{Y} \\
		1 & Y_{n+1}>\overline{Y}.
	\end{cases}
\end{equation*}
If $X_n=1$, we set
\begin{equation*}
	X_{n+1}=\begin{cases} 0 & Y_{n+1}=Y_0 \mbox{ and }\gamma_n<\tau_1\\
		1 & Y_{n+1}>Y_0 \mbox{ and }\gamma_n<\tau_1\\
		2& \gamma_n \ge \tau_1,
	\end{cases}
\end{equation*}
where $\tau_1$ is the stress threshold. Finally, if $X_n=2$ we set
\begin{equation*}
	X_{n+1}=\begin{cases} 0 & Y_{n+1}=Y_0 \mbox{ and }\gamma_n \ge \tau_2\\
		2 & \mbox{otherwise,}
	\end{cases}
\end{equation*}
where $\tau_2$ is the time to repair. Practically, we can assume that, instead of the full information on the stress level $Y_n$, we can only observe a LED $X_n$ with three possible states:
\begin{itemize}
	\item the state $0$ (for instance, the LED is turned off), which tells us that the system is working without any breakdown risk;
	\item the state $1$ (for instance, the LED is turned on), which tells us that the system could break down unless we take the stress level to its default state;
	\item the state $2$ (for instance, the LED is blinking), which tells us that the system is currently broken.
\end{itemize} 
It is also clear, by construction, that the process $X_{n}$ is semi-Markov. Thanks to $X_n$, we can also define the rewards as
\begin{equation*}
	r_n(x,t,a)=\begin{cases} c_1 & a=a_1 \mbox{ and }x \not = 2\\
		c_2 & a=a_2 \mbox{ and }x \not = 2\\
		c_3 & x=2.
	\end{cases}
\end{equation*}
Our goal is to maximize the expected profit obtained up to a fixed horizon $N$ (for instance, if we know that at time $t_N=N$ we must shut down the system for some planned maintenance).\\
If we set $c_1=0$, we can assume that the \textit{slow} regime $a_1$ is actually a \textit{preventive maintenance regime}, in which we shut down the machine to make it cool down. Preventive maintenance represents a quite natural application for semi-Markov decision processes. Indeed, there is no real reason to assume both the time to repair and the system failure time to be exponentially distributed (see \cite{migawa2021semi,sanchez2022optimizing,tomasevicz2006preventive} and references therein). Such a problem is classically approached via a \textit{planning strategy}. Precisely, the optimal time after which the system has to be set in preventive maintenance is calculated a priori. This is equivalent to the application of a certain decision as soon as the system starts working (that is a \textit{change of state}), which cannot be \textit{adjusted} while it is working. Such a practice is perfectly justified by the fact that in these models we do not have any information on the \textit{stress on the system}. On the other hand, even if we have a partial information on it, such a strategy seems to be inaccurate. Clearly, if I can observe, even partially, the \textit{stress level} of the system I should be able to modulate my decision depending on this information. Our sojourn-based approach permits this: we can decide whether to start the preventive maintenance or not depending on the stress level of the system. On the other hand, while the \textit{planning strategy} allows to consider a continuous-time model, the observations of the system cannot be continuous in time, thus we have to resort to discrete-time semi-Markov processes, as described in this example. 
\subsubsection{Numerical results}\label{sec522}
Again, let us consider a specific choice of parameters to exhibit some features of the method. However, before doing this, we need to set the involved distributions. First of all, let us assume that $\tau_1$ is a Poisson random variable with parameter $\lambda$, i.e.
\begin{equation*}
	\bP(\tau_1=s)=\frac{\lambda^s e^{-\lambda}}{s!}, \ s \in \N_0.
\end{equation*}
This is only a simplifying assumption, which is useful to adjust the mode of the distribution according to $\lambda$. A more realistic distribution for $\tau_1$ could be a discrete Weibull one, as introduced in \cite{nakagawa1975discrete}. In both cases, however, $\bP(\tau_1=0)\not = 0$, thus we are admitting possible \textit{instantaneous breakdowns}. Concerning $\tau_2$, we suppose it is a negative binomial random variable of parameters $r \in \N$ and $p \in (0,1)$, i.e.
\begin{equation*}
	\bP(\tau_2=s)=\binom{s+r-1}{r-1}(1-p)^sp^r, \ s \in \N_0.
\end{equation*}
The choice of such a distribution is justified as follows: we can assume the system is composed of $r$ parts that break down simultaneously and that each one of such parts is repaired in a time described according to a geometric distribution of parameter $p$.\\
Concerning $(Y_n)_{n \in \N}$, we set $Y_0=0$ and we suppose that $\bP_{\pi}(Y_{n+1}=y^\prime|Y_n=y)=0$ whenever $|y^\prime-y|>1$, for any policy $\pi \in F^N$, where we omit the initial data in the notation as it is assigned. For fixed $n$, let us distinguish the following sets $\mathcal{E}_n^i \subset \N_0 \times \widetilde{E}$ (where $E=\{0,1,2\}$):
\begin{align*}
	\mathcal{E}_n^1&:=\{(y,x,t): \ \pi_n(x,t)=a_1, \ x \not = 2, y>0\}\\
	\mathcal{E}_n^2&:=\{(y,x,t): \ \pi_n(x,t)=a_1, \ x \not = 2, y=0\}\\
	\mathcal{E}_n^3&:=\{(y,x,t): \ \pi_n(x,t)=a_2, \ x \not = 2\}\\
	\mathcal{E}_n^4&:=\{(y,x,t): x = 2, \ y>0\}\\
	\mathcal{E}_n^5&:=\{(y,x,t): \ x=2, \ y=0\}.
\end{align*}
Clearly $\bigcup_{i=1}^{5}\mathcal{E}_n^i=\N_0 \times \N_0 \times \widetilde{E}$. Thanks to such a decomposition, we can define
\begin{equation*}
	\bP_{\pi}(Y_{n+1}=y+1|Y_n=y,X_n=x,\gamma_n=t)=\begin{cases}
		p_{+,1} & (y,x,t) \in \mathcal{E}_n^1 \cup \mathcal{E}_n^2\\
		p_{+,2} & (y,x,t) \in \mathcal{E}_n^3 \\
		0 & (y,x,t) \in \mathcal{E}_n^4 \cup \mathcal{E}_n^5,
	\end{cases} 
\end{equation*} 
\begin{equation*}
	\bP_{\pi}(Y_{n+1}=y-1|Y_n=y,X_n,\gamma_n)=\begin{cases}
		p_- & (y,x,t) \in \mathcal{E}_n^1\\
		0 & (y,x,t) \in \mathcal{E}_n^2 \cup \mathcal{E}_n^3 \cup \mathcal{E}_n^5\\
		\overline{p}_- & (y,x,t) \in \mathcal{E}_n^4
	\end{cases} 
\end{equation*}
and
\begin{equation*}
	\bP_{\pi}(Y_{n+1}=y|Y_n=y,X_n,\gamma_n)=\begin{cases}
		p_1 & (y,x,t) \in \mathcal{E}_n^1\\
		p_1+p_- & (y,x,t) \in \mathcal{E}_n^2\\
		p_2 & (y,x,t) \in \mathcal{E}_n^3\\
		\overline{p} & (y,x,t) \in \mathcal{E}_n^4\\
		1 & (y,x,t) \in \mathcal{E}_n^5.
	\end{cases} 
\end{equation*}
Now let us focus on the following choice of parameters. We set the horizon $N=200$ and the stress threshold $\overline{Y}=6$. Concerning the distribution of $\tau_1$, we set $\lambda=5$, while we also set $p=0.5$, $r=20$ for the distribution of $\tau_2$. Then we set $p_{+,1}=1/10$, $p_{+,2}=4/5$, $p_-=2/5$, $\overline{p}_-=3/5$, $p_1=1-p_{+,1}-p_-=1/2$, $p_2=1-p_{+,2}=1/5$ and $\overline{p}=1-\overline{p}_-=2/5$. As a first study case, we also consider $c_1=1$, $c_2=5$ and $c_3=-10$.\\
We adopt the same notation of the previous example and we use the $Q$-learning algorithm based on the recursion formula \eqref{Qlearn}. Again, at the $(m+1)$-th iteration and for the $n$-th turn, the action $a$, on which the evaluation $Q_n^{m+1}(x,t,a)$ is taken, is selected in the set $\arg \max Q_{n+1}^m(X_{n+1}^m,\gamma_{n+1}^m,\cdot)$ and $(x,t)=(X_n^m,\gamma_n^m)$. Both $(X_{n}^m,\gamma_{n}^m)$ and $(X_{n+1}^m,\gamma_{n+1}^m)$ are obtained by stochastic simulation, where we also simulate the process $(Y_n^m)_{n \in \N}$ in background. This procedure is repeated for $m=0,\dots,N_{\rm eps}-1$ where $N_{\rm eps}=40000$. Again, for each simulated trajectory $(X_n^m,\gamma_n^m,a_n^m)$ we consider
\begin{equation*}
	R_m=\sum_{n=0}^{N}r(X_n^m,\gamma_n^m,a_n^m)
\end{equation*}
and, to smooth the data, we consider average, minimal and maximal expected reward as
\begin{align*}
	R_{k,{\rm avg}}&=\E\left[\frac{1}{400}\sum_{m=400k}^{400(k+1)-1}R_m\right],\\
	R_{k,{\rm min}}&=\E\left[\min_{400k \le m z  400(k+1)}R_m\right],\\
	R_{k,{\rm max}}&=\E\left[\max_{400k \le m z  400(k+1)}R_m\right].
\end{align*}
The results for different values of the learning rate can be seen in Figure \ref{fig:figure13}. As we anticipated, the average expected reward is increasing with $m$ and, as one can see from the plot, such increase is faster for big values of ${\rm LR}$, as in the previous example. However, this time the gap between the different values of ${\rm LR}$ is much more evident. The cases ${\rm LR}=0.1,0.2$ are evidently slower than the others, while  ${\rm LR}=0.3$, despite being faster than the latter two, performs slightly worse than ${\rm LR}=0.4,0.5$. For this reason, we use as a reference case ${\rm LR}=0.5$, for which the plot of $R_{k, {\rm avg}}$ is shown in Figure \ref{fig:figure14}, together with a confidence interval of $95\%$.
\begin{center}
	\begin{minipage}{0.49\linewidth}
		\includegraphics[width=\linewidth]{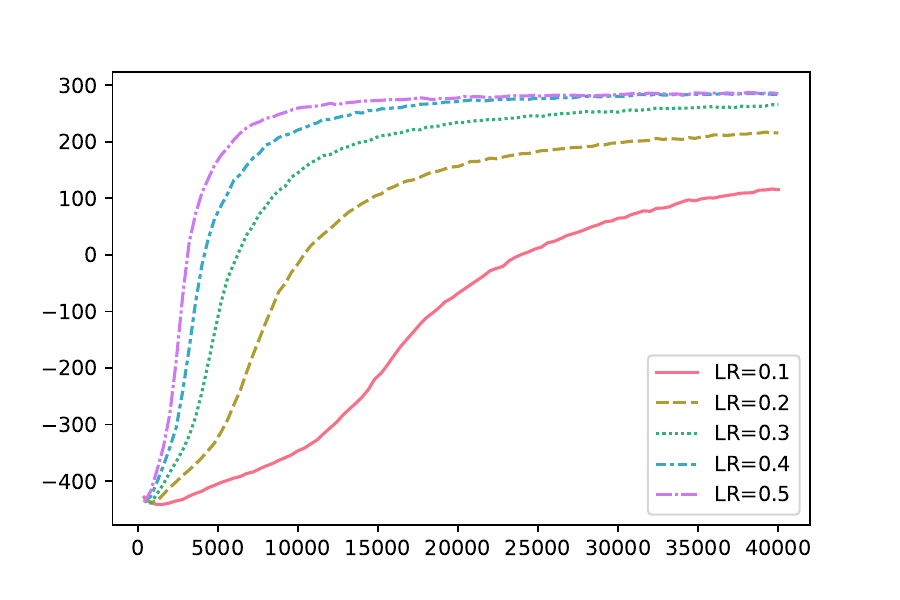}
		\captionof{figure}{Plot of the average expected rewards for different values of the learning rate in the \textit{preventive maintenance} example with $c_1=1$.}
		\label{fig:figure13}
	\end{minipage}
	\begin{minipage}{0.49\linewidth}
		\includegraphics[width=\linewidth]{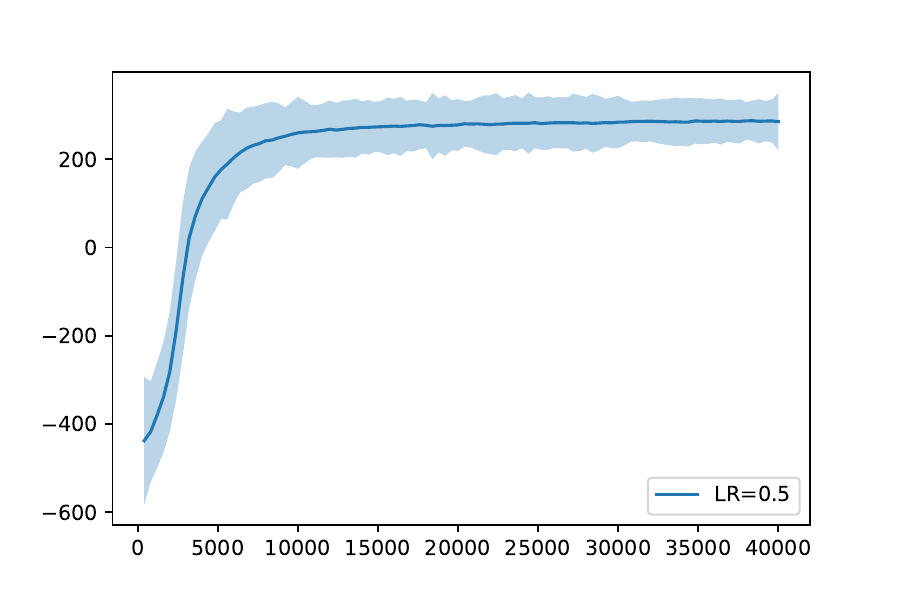}
		\captionof{figure}{Plot of the average expected rewards for ${\rm LR}=0.5$ with a $95\%$ confidence interval in the \textit{preventive maintenance} example with $c_1=1$.}
		\label{fig:figure14}
	\end{minipage}
\end{center}
The plots of the minimal expected rewards for different values of the learning rate are given in Figure \ref{fig:figure15}. Despite the minimal expected rewards are increasing, the gap between the cases ${\rm LR}=0.1, 0.2, 0.3$ and ${\rm LR}=0.4, 0.5$ is again clearly visible with this smoothing procedure. For this reason we consider ${\rm LR}=0.5$ as a reference case, as exploited in Figure \ref{fig:figure16}.
\begin{center}
	\begin{minipage}{0.49\linewidth}
		\includegraphics[width=\linewidth]{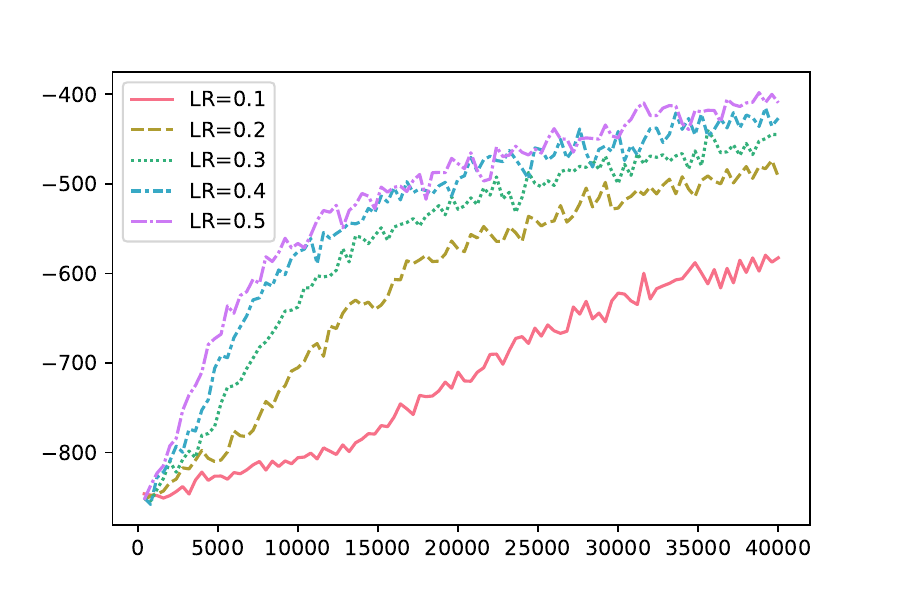}
		\captionof{figure}{Plot of the minimal expected rewards for different values of the learning rate in the \textit{preventive maintenance} example with $c_1=1$.}
		\label{fig:figure15}
	\end{minipage}
	\begin{minipage}{0.49\linewidth}
		\includegraphics[width=\linewidth]{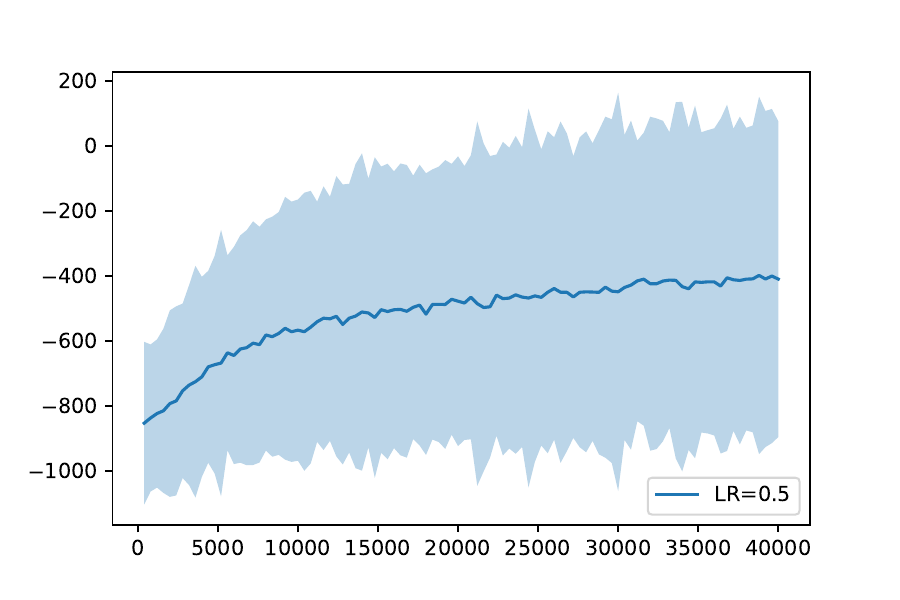}
		\captionof{figure}{Plot of the minimal expected rewards for ${\rm LR}=0.5$ with a $95\%$ confidence interval in the \textit{preventive maintenance} example with $c_1=1$.}
		\label{fig:figure16}
	\end{minipage}
\end{center}
The behaviour of the maximal expected rewards is completely different. Indeed, $R_{k,{\rm max}}$ is first increasing and then decreases smoothly towards a limit value, which seems to be bigger for ${\rm LR}=0.1$. Actually, we expect the limit value to be more or less equal, but the convergence to be much more slower in the case ${\rm LR}=0.1$ with respect to the others. Thus, while on the one hand this case gives the lowest average and minimal expected reward at episode $40000$ with respect to the other ones, on the other hand it gives the highest maximal expected reward, thanks to the fact that the convergence is quite slow. In this case, we could consider ${\rm LR}=0.1$ as a reference case. However, let us exploit again the case ${\rm LR}=0.5$ in Figure \ref{fig:figure18}, so that we can compare the confidence intervals.
\begin{center}
	\begin{minipage}{0.49\linewidth}
		\includegraphics[width=\linewidth]{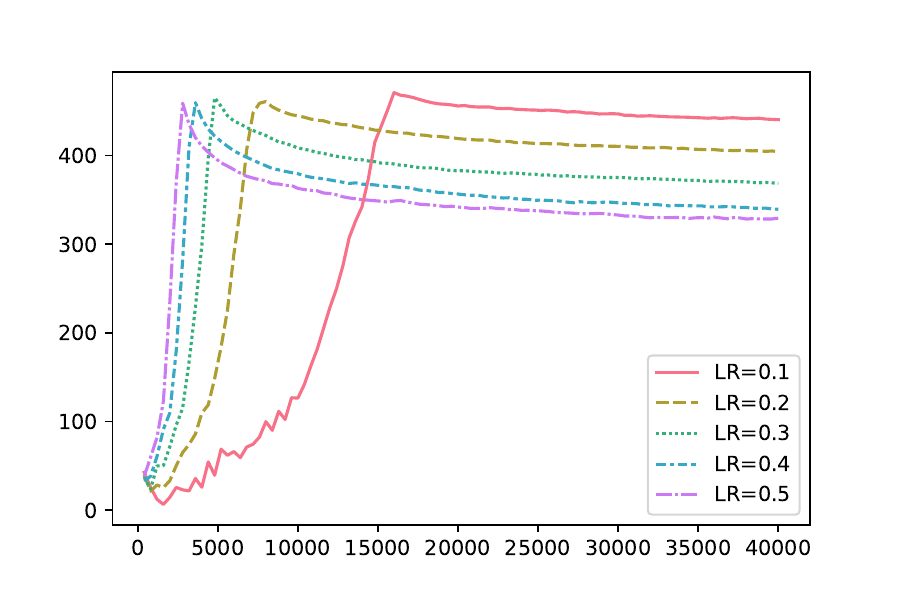}
		\captionof{figure}{Plot of the maximal expected rewards for different values of the learning rate in the \textit{preventive maintenance} example with $c_1=1$.}
		\label{fig:figure17}
	\end{minipage}
	\begin{minipage}{0.49\linewidth}
		\includegraphics[width=\linewidth]{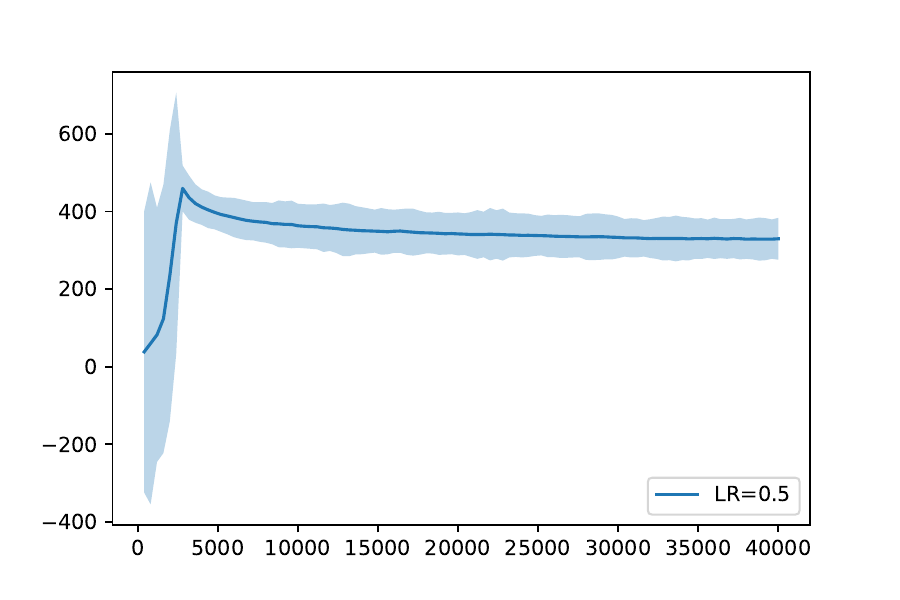}
		\captionof{figure}{Plot of the maximal expected rewards for ${\rm LR}=0.5$ with a $95\%$ confidence interval in the \textit{preventive maintenance} example with $c_1=1$.}
		\label{fig:figure18}
	\end{minipage}
\end{center}
As we already did in the previous example, let us also consider the case of a non-constant learning rate. Precisely, as done in \cite{bhatnagar2021finite}, let us consider
\begin{equation}\label{eq:varLR}
	{\rm LR}(m)=\frac{1}{\lceil \frac{m+1}{5000} \rceil +1},
\end{equation}
which satisfies the \textit{tapering stepsize} assumptions. In Figures \ref{fig:figure19}, \ref{fig:figure20} and \ref{fig:figure21} the plots of the average, minimal and maximal expected rewards respectively are shown for both ${\rm LR}=0.5$ and the learning rate given by \eqref{eq:varLR}. Again, as presumed, the non-constant learning rate case performs slightly worse after a certain number of episodes than the constant one for the average and minimal expected rewards, while in the maximal one it performs slightly better. 
\begin{center}
	\begin{minipage}{0.49\linewidth}
		\includegraphics[width=\linewidth]{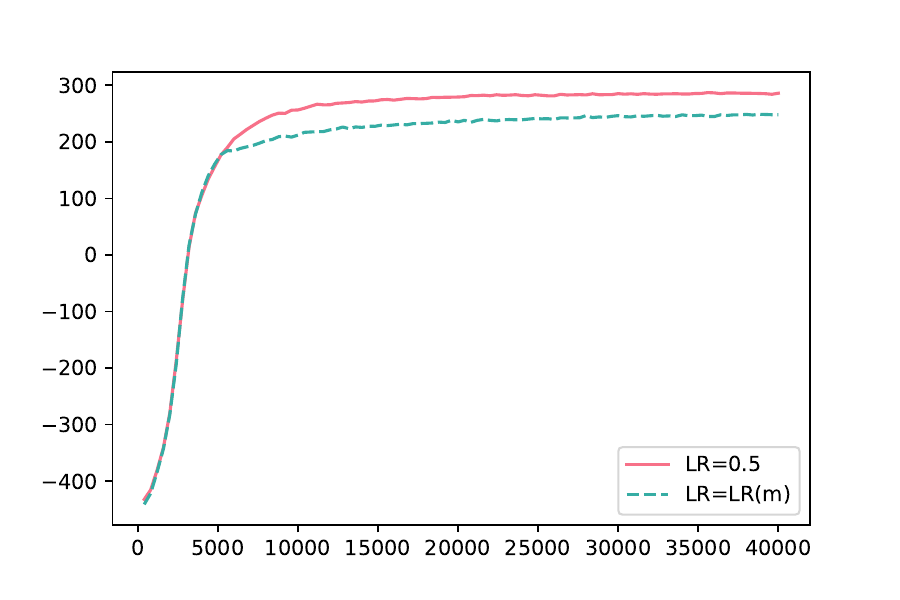}
		\captionof{figure}{Plot of the average expected rewards for fixed and variable learning rates in the \textit{preventive maintenance} example with $c_1=1$.}
		\label{fig:figure19}
	\end{minipage}
	\begin{minipage}{0.49\linewidth}
		\includegraphics[width=\linewidth]{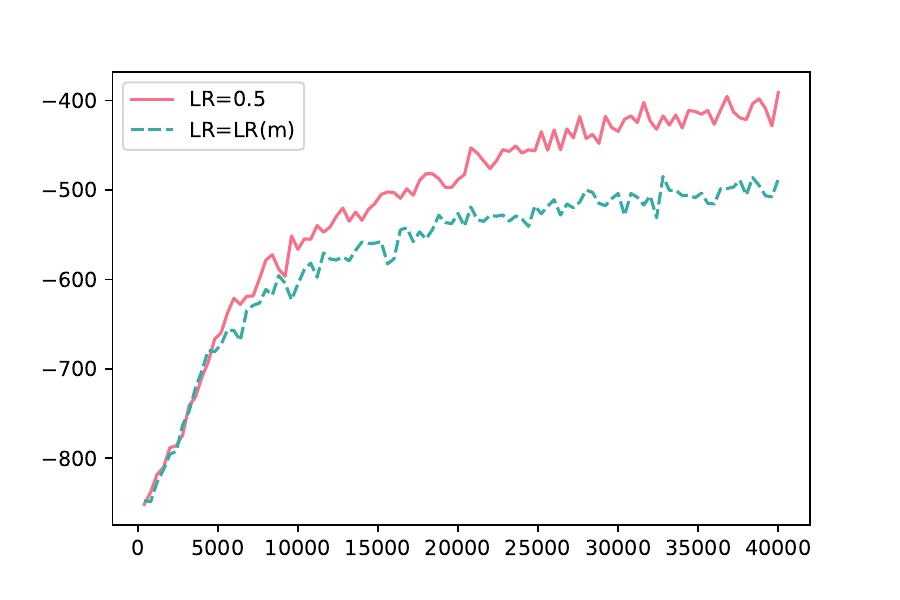}
		\captionof{figure}{Plot of the minimal expected rewards for fixed and variable learning rates in the \textit{preventive maintenance} example with $c_1=1$.}
		\label{fig:figure20}
	\end{minipage}
\end{center}
Let us stress that such a behaviour is connected with the values of $c_1,c_2,c_3$. Indeed, if we consider the \textit{preventive maintenance setting}, i.e. $c_1=0$, while leaving $c_2=5$ and $c_3=-10$, we observe from Figures \ref{fig:figure22}, \ref{fig:figure23} and \ref{fig:figure24} that, despite the behaviour is similar, the choice ${\rm LR}=0.3$ is already sufficient to guarantee a fast and precise approximation. Actually, the case ${\rm LR}=0.3$ performs better than ${\rm LR}=0.5$ in this case: recall that we expect from Proposition \ref{prop:conv1} a \textit{lower} absolute error for smaller values of ${\rm LR}$, while obtaining a possibly slower convergence. Here the price we pay in terms of the speed of convergence is fully balanced by the performance of the algorithm.
\begin{center}
	\begin{minipage}{0.49\linewidth}
		\includegraphics[width=\linewidth]{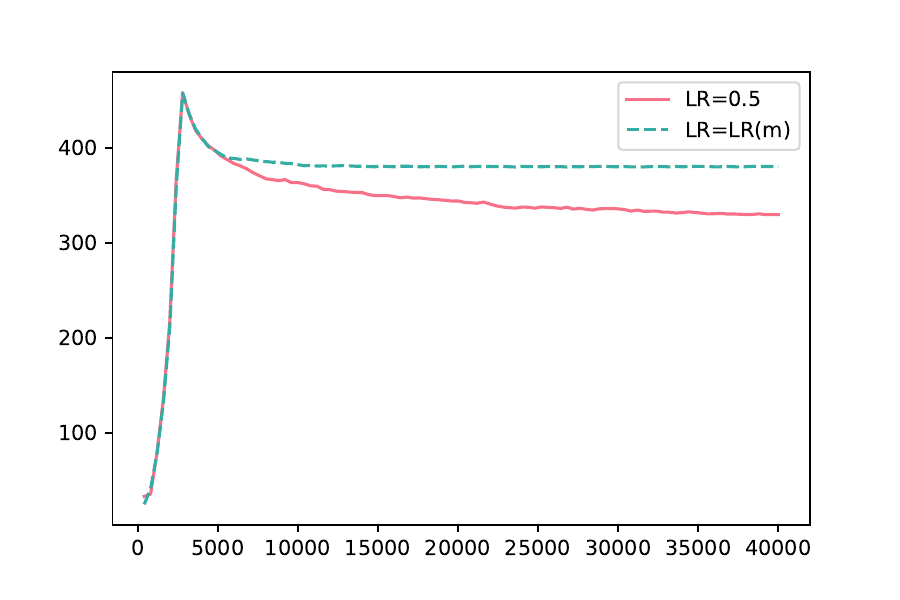}
		\captionof{figure}{Plot of the maximal expected rewards for fixed and variable learning rates in the \textit{preventive maintenance} example with $c_1=1$.}
		\label{fig:figure21}
	\end{minipage}
	\begin{minipage}{0.49\linewidth}
		\includegraphics[width=\linewidth]{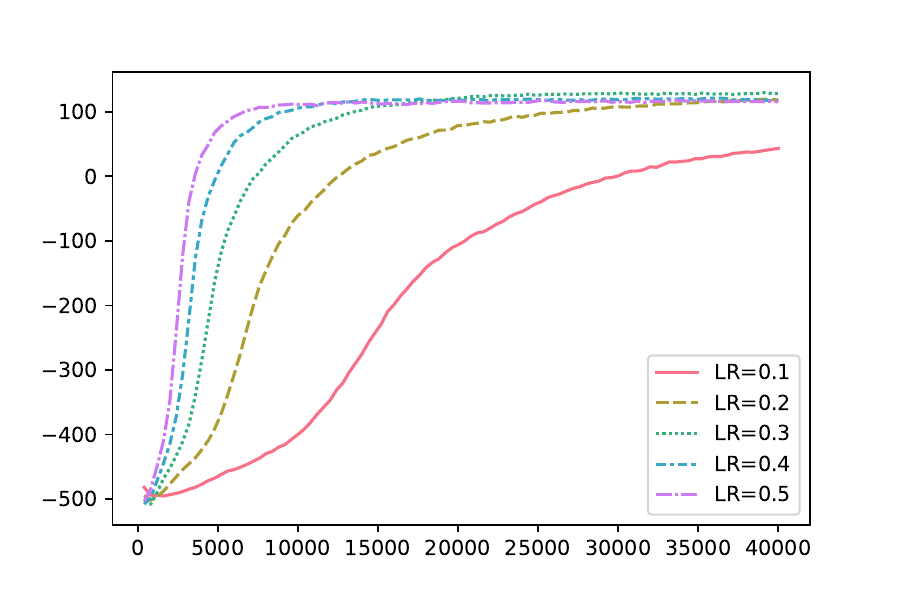}
		\captionof{figure}{Plot of the average expected rewards for different values of the learning rate in the \textit{preventive maintenance} example with $c_1=0$.}
		\label{fig:figure22}
	\end{minipage}
\end{center}
\begin{center}
	\begin{minipage}{0.49\linewidth}
		\includegraphics[width=\linewidth]{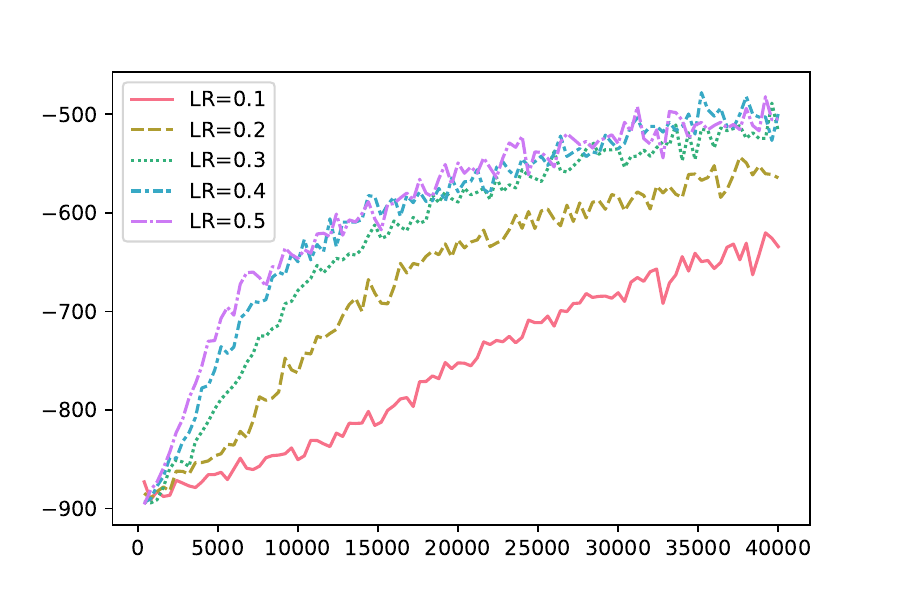}
		\captionof{figure}{Plot of the minimal expected rewards for different values of the learning rate in the \textit{preventive maintenance} example with $c_1=0$.}
		\label{fig:figure23}
	\end{minipage}
	\begin{minipage}{0.49\linewidth}
		\includegraphics[width=\linewidth]{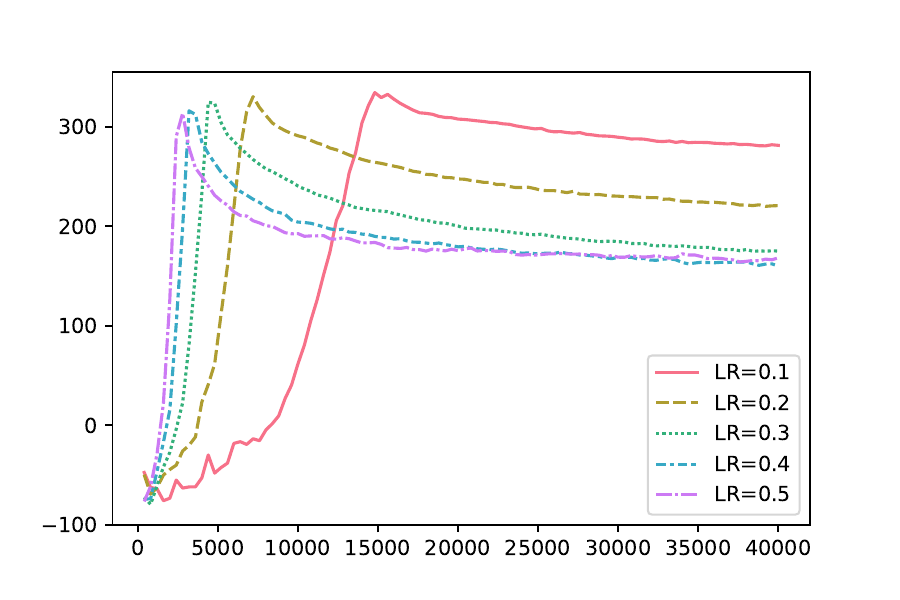}
		\captionof{figure}{Plot of the maximal expected rewards for different values of the learning rate in the \textit{preventive maintenance} example with $c_1=0$.}
		\label{fig:figure24}
	\end{minipage}
\end{center}
Clearly, this also means that the learning rate given by the formula \ref{eq:varLR} should perform better that ${\rm LR}=0.5$ with this choice of parameters, as it should be fast in the first episodes and then lead to almost surely convergence of the method. This is underlined by Figures \ref{fig:figure25}, \ref{fig:figure26} and \ref{fig:figure27}
\begin{center}
	\begin{minipage}{0.49\linewidth}
		\includegraphics[width=\linewidth]{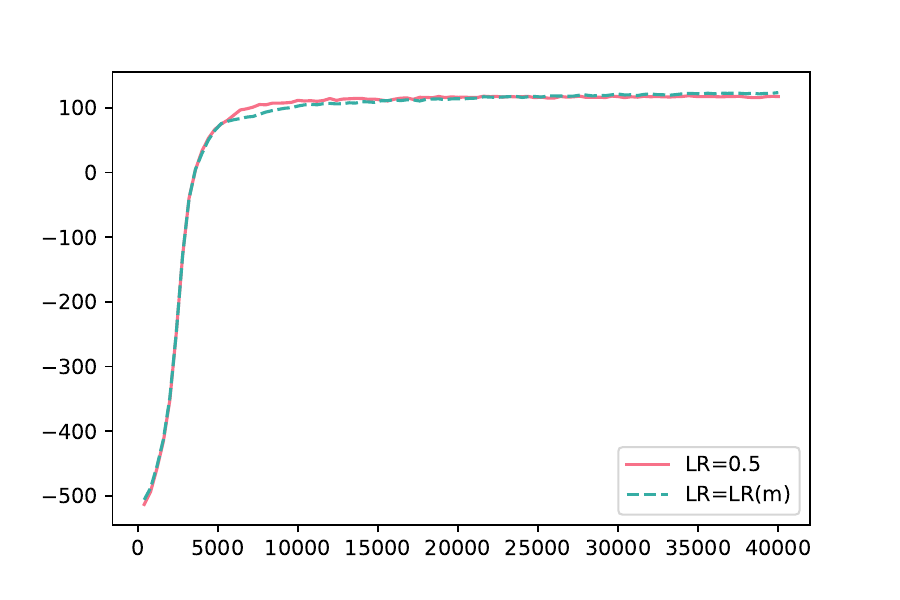}
		\captionof{figure}{Plot of the average expected rewards for fixed and variable learning rates in the \textit{preventive maintenance} example with $c_1=0$.}
		\label{fig:figure25}
	\end{minipage}
	\begin{minipage}{0.49\linewidth}
		\includegraphics[width=\linewidth]{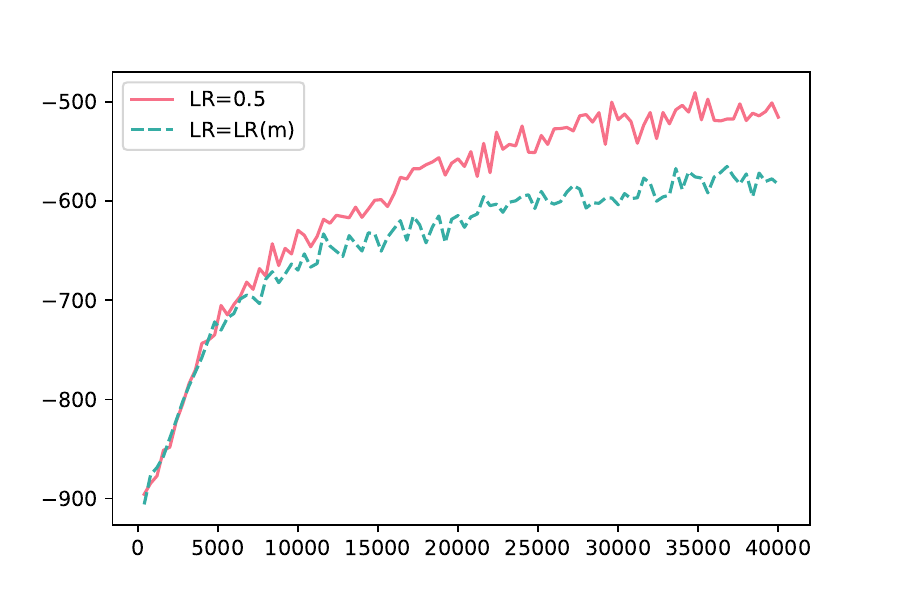}
		\captionof{figure}{Plot of the minimal expected rewards for fixed and variable learning rates in the \textit{preventive maintenance} example with $c_1=0$.}
		\label{fig:figure26}
	\end{minipage}
\end{center}
\begin{center}
	\begin{minipage}{0.49\linewidth}
		\includegraphics[width=\linewidth]{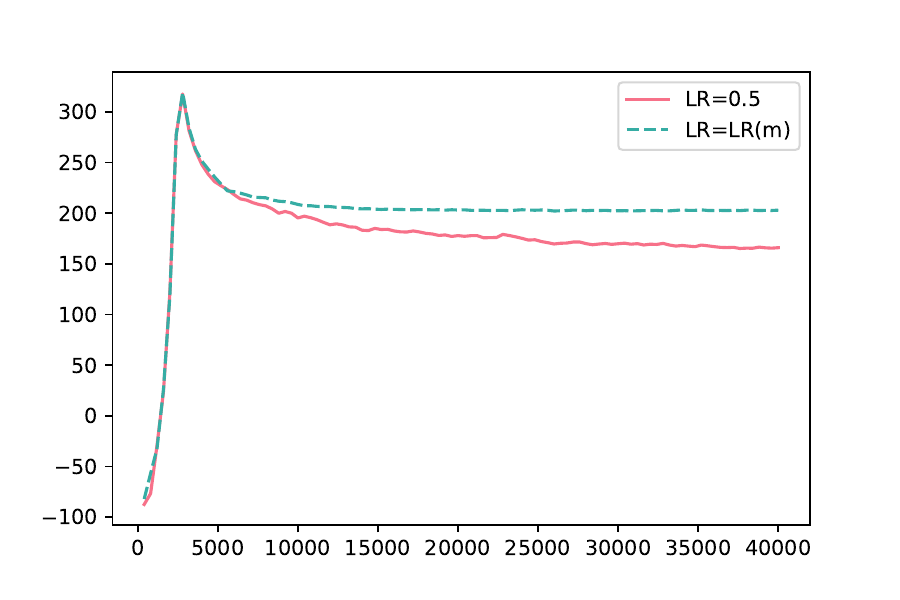}
		\captionof{figure}{Plot of the maximal expected rewards for fixed and variable learning rates in the \textit{preventive maintenance} example with $c_1=0$.}
		\label{fig:figure27}
	\end{minipage}
\end{center}
Let us observe that, in this specific case, convergence is slightly visible at episode $10000$. For this reason, we avoid using the $Q$-network approach in this example.
\section{Discussions and Conclusions}\label{sec6}
In this work we considered an approach to semi-Markov decision processes that is based on the sojourn time process. Precisely, exploiting a characterization of semi-Markov processes in terms of bivariate Markov processes (as in Theorem \ref{thm:char1}), we are able to reduce a semi-Markov decision process to a Markov one and then apply the usual dynamic programming principle. This approach takes in consideration the possibility of observing the system at a discrete set of times that are not necessarily its jump times. Thus, one can decide to change the action dynamically, depending on the present sojourn time. This possibility is apparently not considered in the classical approach to semi-Markov decision processes, as the action is chosen at a jump time and acts as a parameter of the inter-jump time distribution. However, in the discrete setting, we are able to prove another characterization Theorem (precisely, \ref{thm:bij}) that implies the \textit{equivalence} of our approach and the classical one, upon considering a bigger space of actions and suitable parametrizations of the inter-jump times. On the other hand, we think our approach is \textit{more intuitive} whenever the reward function depends directly on the state and the sojourn time at each observation of the process. A biproduct of our first observations is the equivalence between time-changes with Markov additive processes and discrete-time semi-Markov ones. As a consequence, we think that the approach presented in \cite{pachon2021discrete} could be of general interest to determine some generalized Chapman-Kolmogorov equations for the one-dimensional distributions of discrete-time semi-Markov processes. In future works we also want to address the infinite horizon and the optimal stopping cases.\\
The decision problem is not only approached from the theoretical point of view of dynamic programming and Bellman's Equation (in Section \ref{sec3}), but also from the one of the applications (in Section \ref{sec4}). Indeed, in such a context, the transition probability function of the state-sojourn process is usually unknown and different methods have to be implemented in order to determine the optimal policy. While for the infinite horizon setting this is a well established problem, in the finite horizon case this has not been completely addressed. One of the possible strategy is the reduction of episodic tasks to the continuous setting, as exploited in \cite[Section 3.4] {sutton2018reinforcement}, however this is not the case. In \cite{bhatnagar2021finite}, the authors presented a $Q$-learning algorithm adapted to the finite horizon case and proved its convergence properties by using the ODE methods for stochastic recursive relations given in \cite{borkar2000ode}. To ensure the convergence, they used an episode-dependent learning rate. Here we decided to first \textit{sacrifice accuracy for speed}: we used a constant learning rate, for which the convergence cannot be proved, but one can still obtain a bound on the limit mean square error by means the results in \cite{borkar2000ode}. On the other hand, using a constant learning rate permits to study (at least experimentally) the trade-off between accuracy and speed, in analogy of what is proved in \cite[Theorem 2.4]{borkar2000ode}. One can use a first numerical investigation on the constant learning rate case to establish a good form of variable learning rate, so that one gets a suitably fast convergence to the optimal value.\\
In Section \ref{sec5} we tested our approach on two toy examples, that we call the \textit{Switching Coins Problem} and the \textit{Preventive Maintenance with partial observation}. The first problem takes in consideration some natural behaviours of a human interacting with the AI, as we try to consider in the reward function the effect of \textit{gambler's fallacy} (see \cite{warren2018re}). As a consequence, despite the transition probabilities of the underlying model do not depend on the sojourn time process upon choosing an action, the reward function depends on it and thus the decision rules has to be considered as functions of the state and the sojourn time, as in the semi-Markov case. This implies that we cannot exclude that the process obtained via the optimal policy is semi-Markov, despite the underlying model is Markov for constant policies. This example shows that semi-Markov dynamics could appear even in Markov environments if the reward function takes in consideration the sojourn time process, which could be the case to model some human bias.\\
The second problem, instead, is a modification of a classical semi-Markov based one. Indeed, as exploited for instance in \cite{tomasevicz2006preventive}, the fact that the sojourn time in each state is not necessary exponentially distributed tells us that the environment itself is non-Markov, but only semi-Markov. However, in such a model, it is usual to assume that we cannot have any information on the stress to which the system is subject. For this reason, it is clear that we can plan the preventive maintenance \textit{a priori}, i.e. even before starting the system. Here we are considering a situation in which we have some partial information on the fact that the system will break down or not. Such an information can be achieved only through an observation of the system while it is working, thus, in place of planning the preventive maintenance \textit{a priori}, we should decide if we want to send the system in preventive maintenance \textit{in fieri} (i.e. right after the observation). This cannot be done with a classical approach to semi-Markov reinforcement learning, in which the action is taken only in a change of state, but is clearly manageable with our sojourn-based approach. The fact we are considering simple examples permitted us to carry out different tests on the learning rate, so that we could explore their efficiency and then construct an ad-hoc episode-dependent learning rate.\\
Last, but not least, usual $Q$-learning methods require to save all the values of the $Q$ functions for all the states. In our case, this implies that we have to save a certain number of values that grows quadratically with the horizon. As an alternative, one could consider some different function approximation methods. Thus we investigate a naive Deep Q-learning algorithm, that takes not only the state-age process as input but also the actual step. We studied the effect of the variation of some hyperparameters (Batch size, Number of Neurons per dense hidden layer,  Number of Hidden Layers) on the cumulative reward at each episode (actually, on batches of $50$ episodes each). Being this just a first attempt, we clearly obtain a lower reward than the one achieved via the classical $Q$-learning algorithm. However, the results are still promising. In particular, let us stress that the deep $Q$-learning algorithms obtained in literature are usually based on the fact that in the infinite horizon case the function $Q$ is a fixed point of a certain operator and thus it is independent of the step $n$. In the finite horizon case this is no more true and then one should approximate $N$ functions $Q_n$. A first really naive idea could be the implementation of $N$ networks to be trained. This idea is quite unpractical and thus has to be discarded. A second naive idea is to consider $n$ as a further argument of the $Q$ function, i.e. by working directly on the $\mathbf{Q}^*$. This is analogous of working with the strong Markov process $(n,X_n,\gamma_n)$ with a suitable state space $\widetilde{E}_N$ and then defining $Q$ as a function on $\overline{\cD}_N$ (that takes in consideration also the admissible actions). Thus, we introduce here this naive idea to work with finite horizon deep $Q$-learning algorithms. Clearly, a further investigation is needed to improve such an idea. In future works we plan to carry on the study on deep $Q$-learning algorithms in this context (in particular for the finite horizon case, which turned out to have some interesting implications).

\section*{Acknowledgements}
	The authors would like to thank the anonymous referees for their precious comments. The first author is supported by MIUR-PRIN 2017, project Stochastic Models for Complex Systems, no. 2017JFFHSH and by Gruppo Nazionale per l’Analisi Matematica, la Probabilit´a e le loro Applicazioni (GNAMPA-INdAM). The second author is partially supported by INdAM-GNCS, “Research ITalian network on Approximation (RITA)” and UMI Group TAA “Approximation Theory and Applications”.
	
\bibliographystyle{abbrv}
\bibliography{bib}

\end{document}